\documentclass{article}

\usepackage{amssymb,amsmath,makeidx,fancybox,color,colordvi,multicol,graphicx,wrapfig,soul}
\usepackage{cite}

\begin{document}
	
	\begin{center}
		{\large\bf On geometric conditions for reduction of the Moreau sweeping process to the Prandtl-Ishlinskii operator}
	\end{center}
	
	\begin{center}
		D. Rachinskii\footnote{Department of Mathematical Sciences, The University of Texas at Dallas, USA. Email: dmitry.rachinskiy@utdallas.edu}
	\end{center}
	
\begin{abstract}
  The sweeping process was proposed by J. J. Moreau as a general mathematical formalism for quasistatic processes in elastoplastic bodies. This formalism deals with connected Prandtl's elastic-ideal plastic springs, which can form a system with an arbitrarily complex topology. The model describes the complex relationship between stresses and elongations of the springs. On the other hand, the Prandtl-Ishlinskii model assumes a very simple connection of springs.
  This model results in an input-output operator, which has many good mathematical properties. It turns out that the sweeping processes can be reducible to the Prandtl-Ishlinskii operator even if the topology of the system of springs is complex. In this work, we analyze the conditions for such reducibility.	
\end{abstract}

	\section{Introduction}
	Models of quasistatic elasto-plasticity date back to Prandtl's model of an elastic-ideal plastic element, which can be thought of as a cascade connection of an ideal Hook's spring and a Coulomb dry friction element. This simple model accounts for two important effects, saturation of stress with increasing deformation (strain) and hysteresis in the stress-strain relationship. 
	Hysteresis is a manifestation of the fact that stress at a moment $t$ is not a single-valued function of the concurrent deformation value, but rather a function of state of the elasto-plastic material,
	which depends on the history of variations of the deformation prior to the instant $t$.
	Two parameters of Prandtl's model are the stiffness $a$ of the spring and the maximal spring force $r$ (which equals the friction force in the sliding regime for quasistatic deformations).
	
	In order to account for complex relationship between deformation and stress in real materials, Prandtl proposed to model the constitutive law of the material
	with a parallel connection of elastic-ideal plastic elements. A similar idea was developed by Ishlinskii who modeled individual fibers of wire ropes
	by Prandtl's elements. In the Prandtl-Ishlinskii phenomenological model, a finite or infinite set of Prandtl's elements (characterized by different values $(a_i,r_i)$
	of parameters of stiffness and maximal stress) are all subject to the same deformation $\varepsilon(t)$, and the total force (or stress) $\sigma(t)$ is proportional to the sum of all spring forces.
	The operator that maps the time series  $\varepsilon(t)$ of the deformation (input) to 
	the time series of stress  $\sigma(t)$ (output), given a set of initial stresses of all the springs (initial state),
	is known as the Prandtl-Ishlinskii (PI) operator in one-dimensional elasto-plasticity.
	Thanks to the set of good mathematical properties of this operator (see, for example, \cite{mroz,xx2,xx3,xx4,xx5,KP,bs}), 
	its equivalent counterparts have been used in several other disciplines including tribology (Maxwell-slip friction model),
	damage counting and fatigue estimation (the rain flow counting method), and, more recently,
	modeling constitutive laws of smart materials such as piezo-electric and magnetostrictive materials and shape memory alloys.
	One useful property, called the composition rule, is that a composition of PI operators is also a PI operator and, as a consequence, the inverse operator for a PI operator
	is another PI operator. Furthermore, a PI operator and its inverse admit an efficient analytic and numerical implementation.
	This property, in particular, underpins the design of compensation-based control schemes for microactuators and sensors,
	which use smart materials for energy conversion. 
	
	Another fact that facilitates modeling various constitutive laws
	with the PI operator, and is also central to this paper, is stated by the Representation Theorem, which allows one to determine 
	whether a set of input-output data can be modeled by a PI operator and, moreover, equips one with 
	an algorithm for identifying parameters of a PI operator from a simple measurement procedure.
	The Representation Theorem states that if (a) the input-output relationship between deformation and stress 
	is {\em rate-independent}\footnote{Rate-independence means that the operator ${\mathcal P}$ that
		maps the time series $\varepsilon(t)$ of deformation to the time series $\sigma(t)$ of stress commutes with any increasing transformation $\tau(t)$ of time, 
		${\mathcal P} \circ \tau=\tau \circ {\mathcal P}$.
		}; (b) hysteresis loops corresponding to periodic inputs are {\em closed}\footnote{This property is common for most phenomenological 
	models of hysteresis with scalar-valued inputs and outputs including the Preisach model \cite{KP,xx1,xx6} and the Ising model.
	Manifestations of this property are also known as the {\em return point memory}, {\em wiping-out property}, {\em no passing rule},
	and {\em Madelung's update rule}.}; and, (c) every hysteresis loop is {\em centrally symmetric}, then $\sigma(t)=({\mathcal P}_\phi\varepsilon)(t)$,
where ${\mathcal P}_\phi: C(0,T)\to C(0,T)$ is a PI operator.
Furthermore, if properties (a)--(c) are satisfied, then starting from the initial state in which all the springs are relaxed,
applying an increasing input (deformation) $\varepsilon$, and measuring the corresponding increasing output $\sigma=\phi(\varepsilon)$,
one obtains the so called {\em loading curve} $\phi: \mathbb{R}_+\to \mathbb{R}_+$, which completely defines 
the PI operator ${\mathcal P}_\phi: C(0,T)\to C(0,T)$  on the class of all continuously varying inputs $\varepsilon: [0,T]\to\mathbb{R}$
via an explicit formula involving $\phi$ and a sequence of local running extremum values of $\varepsilon$ (see \eqref{triangles} below). 
It is worth mentioning that the same measurement of the loading curve is used in nonlinear elasticity 
to identify the non-hysteretic relationship between $\varepsilon$ and $\sigma$ by the Nemytskii operator $\sigma(t)=\phi(\varepsilon(t))$.

The main premise of the of the phenomenological PI model is that Prandtl's elements do not interact.
The model of Moreau addresses a much more general setting, in which an arbitrary spatial configuration of nodes
connected by Prandtl's elastic-ideal plastic springs deforms quasistatically (the balance of forces at each node is zero at all times)
in response to either (i) external forcing applied at a selected set of nodes or (ii) controlled variation of the distance between nodes for a selected set of pairs of nodes, or 
(iii) simultaneous application of the two types of inputs above. In this setting, the relationship between the vector of stresses of individual springs and the input
is complicated and is described by the differential inclusion (or, equivalently, a variational inequality), which is known as the Moreau {\em sweeping process}.
This process has a nice kinematic interpretation, in which variations of the input  induce 
the motion and deformation of the convex domain $Z=Z(t)$ of admissible stresses,
and the set $Z$ drags a point representing the vector of stresses according to a natural kinematic rule.

Here, we consider an (arbitrary) configuration of Prandtl's elastic-ideal plastic springs
controlled by one scalar-valued input-the distance $g(t)$ between two selected nodes, $A$ and $B$.
The corresponding sweeping process belongs to a special class known as a (multidimensional) play operator
with a unidirectional input. Even for this setting, the relationship between the input $g(t)$ and the stress $\sigma_i(t)$
of any given spring can be quite complicated. In particular, hysteresis loops corresponding to a periodically varying $g$ can be non-closed,
thus violating property (b) of the PI operator stated above.
On the other hand, there are multiple examples of topologies of the spring configuration for which 
the relationship between $g$ and $\sigma_i$ is a PI operator ${\mathcal P}_{\phi_i}$ for every $i$
(as well as the total reaction force at the node A, and at the node B, is related to $g$ via a PI operator).
In particular, such examples can be constructed using the composition property of PI operators.

In this paper, we consider the following question: Under which assumptions the Moreau sweeping process driven 
by a scalar-valued input $g(t)$ is equivalent to a PI operator? The answer will be given in geometric terms involving
a vector-valued analog of the loading curve, which we define for the sweeping process.
We show that if this curve and the domain $Z$ of admissible stresses satisfy simple geometric conditions,
then the reaction force of the system of coupled Prandtl's elements responds 
to the input $g(t)$ exactly in the same way as the reaction force of a system of decoupled Prandtl's elements of a certain effective PI model.
Furthermore, $\sigma_i(t)=({\mathcal P}_{\phi_i} g)(t)$ for each spring of the Moreau model,
where $\phi_i$ is an appropriate projection of the loading curve of the Moreau model
in the space of stresses. In other words, we perform a reduction of the more complex model of Moreau
to a simpler affective PI model when such reduction is possible.

The paper is organized as follows. In the next section, we present an outline of models of Moreau and Prandtl-Ishlinskii.
In Section \ref{main}, the main reduction theorem is proved. Section \ref{discussion} contains some discussion.
In particular, we show that any Moreau model obtained by a small perturbation
of a PI model satisfies the conditions of the main theorem and hence is reducible to to an effective PI model.
The paper ends with a summary of results and a few concluding remarks.

\section{The model of Moreau}

\subsection{Mechanical setting}\label{1.1}
\begin{figure}[h]
	\begin{center}
		\includegraphics[width=9truecm]{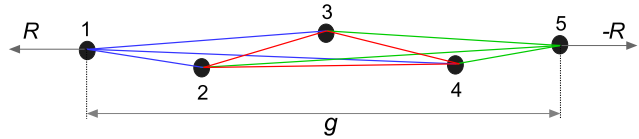}
	\end{center}
	\vskip-10truemm
	%\centerline{(a) \hskip5.5truecm (b)}
	%\vskip-2truemm
	\caption{\small Connection of 5 nodes with 9 elastic-ideal plastic springs.}
	\label{fig1}
\end{figure}

Following the model of Moreau, let us consider an arrangement of $N$ nodes, some of which are connected by elastic-ideal plastic springs, see Fig.~\ref{fig1}.
For simplicity, we consider a one-dimensional arrangement (see Remark \ref{??}). That is,
all the nodes lie on a straight line (the $x$-axis) and the springs, as well as spring forces, are elongated along this line at all times.
The coordinate of node $i$ will be denoted $x_i$.
It is assumed that there is a configuration $(x_1^*,\ldots,x_N^*)$ called the zero configuration, in which the springs experience zero stress,
and we consider deformations of the springs relative to this configuration. The deformation of the spring connecting nodes $i$ and $j$ is denoted by $\varepsilon_{ij}$:
\begin{equation}\label{config}
\varepsilon_{ij}=x_j-x_i-(x_j^*-x_i^*),
\end{equation}
and the force (stress) of this spring is denoted by $\sigma_{ij}$.

This arrangement of springs can be associated with a non-directed graph $\Gamma$ with the nodes associated with graph's vertices and springs associated with graph's edges.
Without loss of generality, this graph is assumed to be connected.
By $M$, we denote the set of indices $(ij)$ such that $(ij)\in M$ whenever there is
a spring connecting the nodes $i$ and $j$. Here and henceforth, we agree that $(ij)=(ji)$
and $i\ne j$ for $(ij)\in M$.

The general Moreau sweeping process models a quasistatic response of the system of springs to external controls (inputs).
There are two types of admissible controls: (a) an external force $F_i=F_i(t)$ applied to a node $i$;
and, (b) a {\em moving affine constraint} that prescribes a distance $L_{ij}=x_j^*-x_i^*+g_{ij}(t)$ between a pair of nodes $i$ and $j$ at all times.
Multiple moving constraints and external forces may be applied simultaneously at several pairs of nodes and nodes, respectively.
However, in this paper, we consider a system with one control. To be specific, we choose to consider a system of springs
controlled by one moving constraint. Without loss of generality, we assume that this constraint defines the distance
between nodes $1$ and $N$ and that these nodes are not connected by a spring:
\begin{equation}\label{input}
\varepsilon_{1N}=g(t)
\end{equation}
where $g\in W^{1,1}$ is a given function of $t\ge 0$ satisfying $g(0)=0$ and $(1\,N)\not \in M$.

We note that Moreau's theory equally applies to  configurations,
where a pair of nodes can be connected by multiple springs. In this case,
$\Gamma$ is a multigraph.

\subsection{Prandtl's elastic-ideal plastic spring}\label{1.2}
\begin{figure}[h]
	\begin{center}
		\includegraphics[width=8truecm]{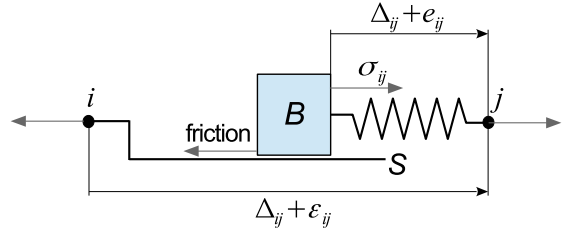}
	\end{center}
	\vskip-10truemm
	%\centerline{(a) \hskip5.5truecm (b)}
	%\vskip-2truemm
	\caption{\small  A nonlinear spring connecting two nodes in the model of Moreau can be represented as a combination
		of a dry friction element and an ideally elastic spring obeying the Hook's law \eqref{hook}. The deformation $e_{ij}$ of the ideal spring is called elastic deformation;
		$\varepsilon_{ij}$ is the elongation of the distance between the nodes $i$ and $j$ (with respect to the distance $\Delta_{ij}=x_j^*-x_i^*$ at zero configuration).
		Dry friction between the box $B$ and the surface $S$ produces the friction force, which is opposite to the ideal spring force at all times (the quasistatic model). The magnitude of the friction force is limited by the maximal value $a_{ij}\rho_{ij}$. Therefore, the ideal spring deforms but the box does not move with respect to the surface $S$ as long as $|\sigma_{ij}|<a_{ij}\rho_{ij}$.
		When $|\sigma_{ij}|=a_{ij}\rho_{ij}$, the deformation $e_{ij}$ of the ideal spring and the force remain constant, while the box moves with respect to the surface in the direction
		of the spring force. }
	\label{fig2}
\end{figure}

In the models of Moreau and Prandtl-Ishlinskii, the stress $\sigma_{ij}=\sigma_{ij}(t)$ and deformation $\varepsilon_{ij}=\varepsilon_{ij}(t)$ of each spring are related by
Prandtl's nonlinear hysteretic constitutive law that combines an ideally elastic spring with a dry friction element, see Fig.~\ref{fig2}.
According to this constitutive law,
\begin{equation}\label{hook}
\sigma_{ij}=a_{ij} e_{ij},
\end{equation}
where
$a_{ij}>0$ is the Young's modulus of the spring in the elastic domain, and the internal variable $e_{ij}$ called {\em elastic deformation}
is related to the deformation $\varepsilon_{ij}$ by the so-called one-dimensional {\em stop} operator $\mathcal{S}_{\rho_{ij}}$:
\begin{equation}\label{stop}
%\sigma_{ij}=a_{ij} e_{ij},\qquad
e_{ij}=\mathcal{S}_{\rho_{ij}}[\varepsilon_{ij}].
\end{equation}
Here $\rho_{ij}$ is the maximal value of the elastic deformation, which defines the maximal magnitude $r_{ij}:=a_{ij} \rho_{ij}$ of stress for the spring,
{\em i.e.}~%Hence, the stress is limited to the  interval
\begin{equation}\label{adm}
|\sigma_{ij}|\le r_{ij}
\end{equation}
at all times.

\begin{figure}[h]
	\begin{center}
		\includegraphics[width=8truecm]{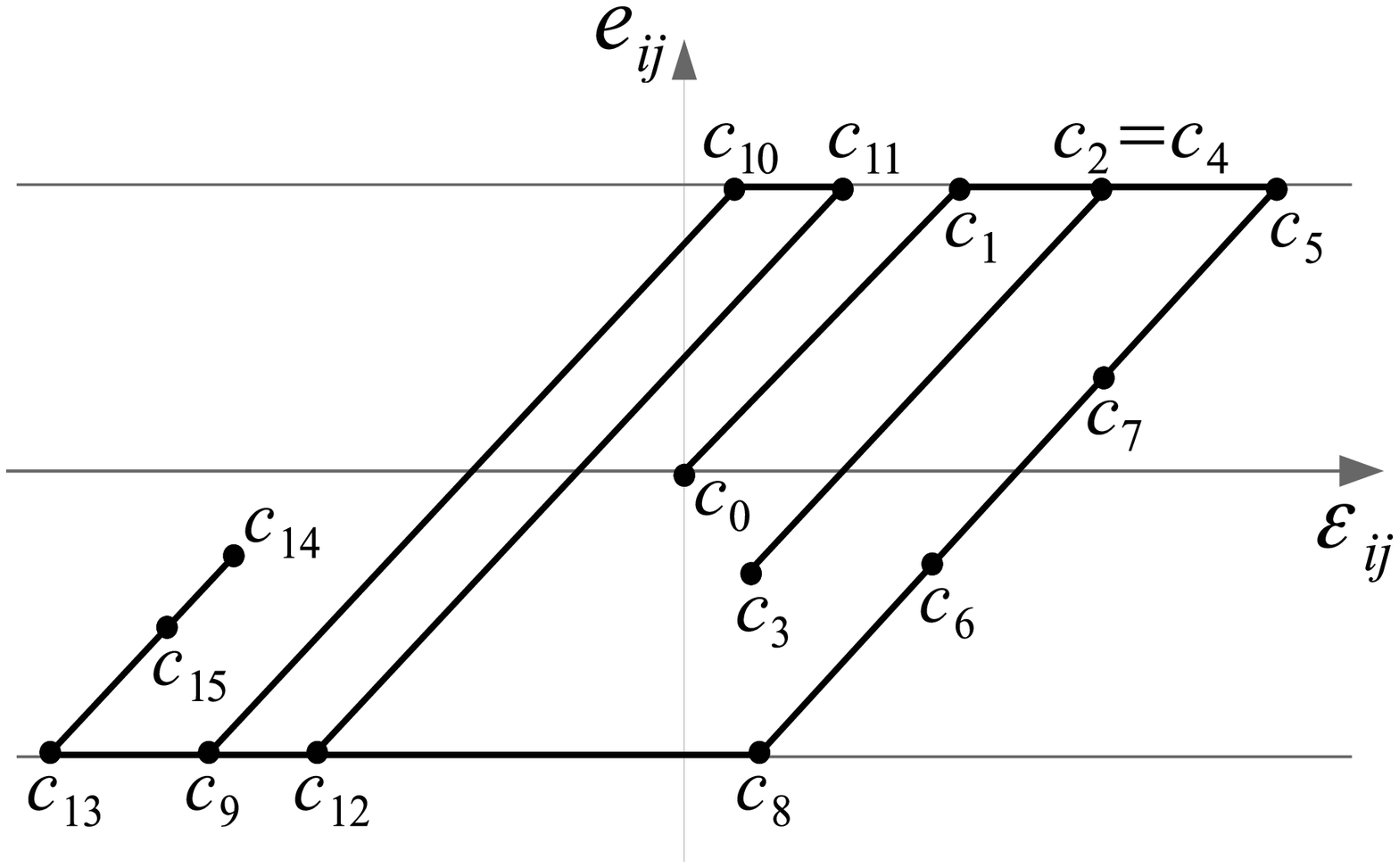}
	\end{center}
	\vskip-10truemm
	%\centerline{(a) \hskip5.5truecm (b)}
	%\vskip-2truemm
	\caption{\small
		An example of an input-output trajectory $c_0\,c_1\,\ldots\,c_{15}$ of a nonlinear spring
		with $c_n=(\varepsilon_{ij}(t_n),e_{ij}(t_n))$ and $t_0<t_1<\cdots$. The time series of the deformation $\varepsilon_{ij}$
		and the elastic deformation $e_{ij}$ are related by the stop operator $e_{ij}={\mathcal S}_{\rho_{ij}}[\varepsilon_{ij}]$.
		Slanted lines have the slope $1$; horizontal lines are $e_{ij}=\pm \rho_{ij}$.
	}
	\label{fig3}
\end{figure}

For piecewise monotone inputs $\varepsilon_{ij}$, relationship \eqref{stop} between the time series of the deformation and the elastic deformation %of a spring
is given by the explicit formula%s
\[
%\sigma_{ij}(t)=a_{ij} e_{ij}(t), \qquad
e_{ij}(t)=F_{\rho_{ij}}\left(\sigma_{ij}(t_n)+\varepsilon_{ij}(t)-\varepsilon_{ij}(t_{n})\right),\quad \ t\in[t_n,t_{n+1}],
\]
where $0=t_0<t_1<t_2<\cdots$ is a partition of the time domain into intervals of monotonicity of the input $\varepsilon_{ij}$
and
\[
F_\rho(x)=\left\{
\begin{array}{rc}
-\rho,& x<-\rho,\\
x,& -\rho\le x\le \rho,\\
\rho,& x>\rho,
\end{array}
\right.
\]
see Fig.~\ref{fig3}. Furthermore, as shown in \cite{KP}, operator \eqref{stop} admits a continuous extension
from the dense set of piecewise monotone inputs to the whole space in each of the spaces $C$, $W^{1,1}$, and $BV$.
In the space $W^{1,1}$, operator \eqref{stop} can be alternatively defined by the variational inequalities 
\begin{equation}\label{var}
|e_{ij}(t)|\le \rho_{ij};\qquad  (\dot \varepsilon _{ij}(t) - \dot e_{ij} (t)) (e_{ij}(t)-z)\ge 0 \ \ {\rm for\ all} \ \ |z|\le \rho_{ij},
\end{equation}
which should be satisfied for a.e.~$t$. This law postulates the ideal elastic response of the spring as long as the stress remains in the range $-r_{ij}<\sigma_{ij}<r_{ij}$.
If the maximal admissible value of the stress is reached but the deformation continues to increase in absolute value,
the stress remains at its maximal value, {\em i.e.} $\sigma_{ij}=r_{ij}$ for $\dot \varepsilon _{ij}>0$ and $\sigma_{ij}=-r_{ij}$ for $\dot \varepsilon _{ij}<0$,
that is the spring {\em yields}.

Everywhere below, we assume that all the deformations and forces are initially zero, that is $\varepsilon_{ij}=e_{ij}=0$, $\sigma_{ij}=0$ for all $(ij)\in M$ at the initial moment $t_0=0$.

%\vskip-10truemm

\subsection{Configuration space}
Define the vector $\varepsilon$ of deformations $\varepsilon_{ij}$ (where $(ij)\in M$) of all the springs
and the vector $\sigma$ of their corresponding stresses $\sigma_{ij}$.
Similarly, $e$ denotes the vector of elastic deformations $e_{ij}$, and we will also
use the vector
\[
p:=\varepsilon-e
\]
of {\em plastic deformations} $p_{ij}$.
The model of Moreau defines dynamics in the configuration space $\mathbb{ E}=\mathbb{ R}^m\ni \varepsilon, \sigma, e, p$, where $m=\# M$ is the total number of springs.

In what follows, the set of geometric constraints \eqref{config} with the additional affine moving constraint \eqref{input} will be expressed as the inclusion
\begin{equation}\label{eps}
\varepsilon (t) \in W +k_0 g(t),\qquad t\ge0,
\end{equation}
with an appropriate choice of the subspace $W$ of the configuration space $\mathbb{E}$
and a vector $k_0$.
%where $W$ is a subspace of the configuration space $\mathbb{E}$ and, without loss of generality,
%we can assume that the vector $k_0$ is orthogonal to $W$ with respect to the standard scalar product
%\[
%\langle y, z\rangle =\sum_{(ij)\in M} y_{ij}z_{ij},\qquad y, z\in\mathbb{E}.
%\]
%We remark that the set of geometric constraints \eqref{config} defines the subspace ${\rm span}\, \{W,k_0\}$ of $\mathbb{E}$ parametrically.
We remark that geometric constraints in parametric form \eqref{config}
can be equivalently expressed in terms of the cycles $(i_1 i_2 \cdots i_k)$ of the graph $\Gamma$
with the added edge $(1\,N)$. Namely,
the sum of deformations along the edges $(i_1i_2), (i_2i_3),\ldots, (i_ki_1)$ of every cycle must be zero:
\begin{equation}\label{eps'}
\varepsilon_{i_1 i_2} +\varepsilon_{i_2 i_3}+\cdots +\varepsilon_{i_{k-1}i_k}+\varepsilon_{i_k i_1}=0.
\end{equation}

The main assumption of Moreau model is that deformations are quasistatic.
Therefore, the balance of forces at each node is zero at any moment:
\begin{equation}\label{sig}
\sum_{(ij)\in M} \sigma_{ij} =0 \quad {\rm for\ every} \quad i=2,\ldots,N-1,
\end{equation}
where the summation is over $j$. The balance of forces at each of the nodes $1$ and $N$
connected by the moving constraint includes the (unknown) reaction $R$ of this constraint:
\begin{equation}\label{sig''}
R+\sum_{(1j)\in M} \sigma_{1j} =-R +\sum_{(Nj)\in M} \sigma_{Nj} =0.
\end{equation}
Relations \eqref{sig''} follow from \eqref{sig}.

\subsection{Moreau sweeping process}
Let us briefly summarize the results of Moreau for the particular case
of a system of springs with one control \eqref{input} that we are considering.
For a closed convex set $Z\subset \mathbb E$, we denote by
$N_Z(z)$ the external normal cone to the set $Z$ at a point $z\in Z$:
\begin{equation}\label{e5}
N_Z(z)=\{y\in \mathbb{E}:\ \langle y,z-\tilde z \rangle \ge 0 \ \ {\rm for\ all} \ \ \tilde z\in Z\},
\end{equation}
where the standard scalar product
\begin{equation}\label{sp}
\langle y, z\rangle =\sum_{(ij)\in M} y_{ij}z_{ij},\qquad y, z\in\mathbb{E},
\end{equation}
is used.
Let us introduce the parallelepiped
$$
C=\{\sigma\in \mathbb{E}:\ |\sigma_{ij}|\le r_{ij}, \ (ij)\in M\}
$$
of admissible stresses defined by \eqref{adm}.
As was shown by J. J. Moreau, system
%Moreau shows that the evolution of the system of springs at all times is described by the relationships
\eqref{config}--\eqref{stop}, \eqref{sig} is equivalent to the system
\begin{eqnarray}
&&e+p\in W + k_0 g(t), \label{e1}\\
&&\sigma=Ae, \label{e15}\\
&&\sigma\in W^\perp \cap C, \label{e2}\\
&&\dot p \in N_C(\sigma) \label{e3}
\end{eqnarray}
with time dependent variables $e,p,\sigma\in \mathbb{E}=\mathbb{R}^m$.
Here $A$ is the diagonal positive matrix with the entries $a_{ij}$ such that \eqref{e15} is equivalent to \eqref{hook};
recall that $e+p=\varepsilon$, hence \eqref{e1} is nothing else as the relationship \eqref{eps},
which is equivalent to the set of geometric constraints \eqref{config} with the additional
moving affine constraing \eqref{input}; the inclusion $\sigma\in W^\perp$, where $W^\perp$ denotes the orthogonal complement
of $W$ in the space $\mathbb{E}$, is equivalent to equations \eqref{sig} that state the balance of forces;
and, the differential inclusion $\dot p \in N_C(\sigma)$, where
dot denotes differentiation with respect to time and $\sigma\in C$, is equivalent to the
variational inequality \eqref{var} that expresses the constitutive law \eqref{stop}
of Prandtl's spring.

Following J. J. Moreau, it is convenient to use the rescaled variables
$$
u=A^\frac12 e,\qquad v=A^\frac12 p,
$$
where the diagonal matrix $A^\frac12$ is the positive square root of $A$.
Introducing the orthogonal subspaces
\begin{equation}\label{UV}
U:= A^\frac12 W,\qquad V:=A^{-\frac12}W^\perp
\end{equation}
and the scaled parallelepiped
\begin{equation}\label{pi}
\Pi:=A^{-\frac12} C,
\end{equation}
and noticing that
$$
\langle \dot p, \sigma-\tilde z\rangle=\langle\dot p, Ae-\tilde z\rangle=\langle A^\frac12 \dot p, A^\frac12 e-A^{-\frac12}\tilde z\rangle=\langle \dot v, u-A^{-\frac12}\tilde z\rangle,
$$
we see that equations \eqref{e1}--\eqref{e3} are equivalent to the system
\[
\begin{array}{rcl}
&&u + v\in U +  A^\frac12 k_0 g(t),\\ %\label{f1}\\
&&u \in V \cap \Pi,\\ % \label{f2}\\
&&\dot v \in N_\Pi (u). %\label{f3}
\end{array}
\]
%Note that $V=U^\perp$.
Now, introducing the orthogonal projection
\begin{equation}\label{f0}
f_0={\rm proj}_V\, A^\frac12 k_0
\end{equation}
of the vector $A^\frac12 k_0$ on the subspace $V=U^\perp$, %(along the subspace $U=V^\perp$),
%we can rewrite Eq.~\eqref{f1} equivalently as
%$u + v \in U +  f_0 g(t)$.
and using the new variables
$$
s = u- f_0 g(t), \qquad \xi = u+v-f_0 g(t),
$$
one %obtain from system \eqref{f1}--\eqref{f3} the system
can rewrite this system as
\begin{eqnarray}
&&\xi\in U, \label{g1}\\
&&s \in V \cap \Pi -f_0 g(t) \label{g2},\\
&&\dot \xi - \dot s \in N_{\Pi - f_0 g(t)} (s). \label{g3}
\end{eqnarray}
Finally, combining inclusion \eqref{g3} with
$
-\dot \xi \in U =V^\perp=N_V(s)
$
and using the identity $N_{D_1}(s)+N_{D_2}(s)=N_{D_1\cap D_2}(s)$, we arrive at the differential inclusion
\begin{equation}\label{g5}
- \dot s \in N_{\Pi \cap V - f_0 g(t)} (s),
\end{equation}
which is known as the {\em Moreau sweeping process} with the characteristic set (input) $Z(t)=\Pi \cap V - f_0 g(t)$.

A few remarks are in order.
%We note that
First, $\Pi \cap V$ is a centrally symmetric convex polytope.
Its central symmetry is important for the following discussion.
Second, the characteristic set $Z(t)=\Pi \cap V - f_0 g(t)$ of the Moreau process is obtained as time-dependent shift of this polytope.
The Moreau process with a characteristic set of the form $Z(t)=Z_0+y(t)$, where $y$ is a single-valued function
$y:\mathbb{R}_+\to V$ and $Z_0\subset V$ is a convex set,
is known as the multi-dimensional play operator with input $y$ \cite{KP}.
In our case, the input $y(t)=-f_0 g(t)$ has a fixed direction $f_0\in V$, {\em i.e.}
we consider the Moreau process of the play type, which in effect has a one-dimensional input.
In general, for a system with multiple inputs (moving constraints and external forces),
the shape of the characteristic set $Z(t)$ may change with time, {\em i.e.}
$Z:\mathbb{R}_+\to 2^V$ is a more general regular set-valued convex-valued function.

Eq.~\eqref{g5} is equivalent to the differential inclusion
\begin{equation}\label{s}
- \dot u + f_0 \dot g(t) \in N_{\Pi \cap V} (u)
\end{equation}
known as the multi-dimensional stop operator with the input $f_0 g(t)$.
This differential inclusion coupled with the initial condition $u(0)=0$
has a unique solution $u\in W^{1,1}(0,T; V)$ for any input
$g\in W^{1,1}(0,T; \mathbb{R})$; regularity properties of the solution operator
that maps $g$ to $u$ in spaces $W^{1,1}$, $C$, and $BV$ are well understood.
Furthermore, system \eqref{g1} -- \eqref{g3} coupled with the initial conditions $s(0)=\xi(0)=0$
(we also assume $g(0)=0$) has a solution $(\xi,s)$. While forces of the springs (the component $s$) are defined uniquely for a given input $g=g(t)$,
simple examples show that the solution $(\xi,s)$ may be non-unique, that is deformations are not necessarily uniquely defined\footnote{The simplest example is a system
	of two identical springs connecting nodes 1 and 2 and nodes 2 and 3, respectively, with the moving constraint applied to nodes 1 and 3.}.
It should be noted that examples with multiple solutions are non-generic. However, the author is not aware of results that would establish uniqueness
of deformations under genericity assumptions.

\subsection{Example}
Consider the system of springs shown in Fig.~\ref{fig1}.
%The elongation of the distance $\varepsilon_{15}=g(t)$ between the nodes 1 and 5 is the input of the system.
%All the processes are assumed to be quasistatic, that is the balance of all forces (including reaction forces due to constraints)
%is zero at every node at all times.
The configuration space $\mathbb{E}=\mathbb{R}^9$ of this system
consists of vectors
$\varepsilon=(\varepsilon_{12},\varepsilon_{13},\varepsilon_{14},\varepsilon_{23},\varepsilon_{24},\varepsilon_{25},\varepsilon_{34},\varepsilon_{35},\varepsilon_{45})'$ (where prime denotes the transposition),
which satisfy the following set of geometric constraints:
$$
\varepsilon_{23}=\varepsilon_{13}-\varepsilon_{12}, \ \varepsilon_{24}=\varepsilon_{14}-\varepsilon_{12}, \ \varepsilon_{25}=\varepsilon_{15}-\varepsilon_{12}, \
\varepsilon_{34}=\varepsilon_{14}-\varepsilon_{13}, \ \varepsilon_{35}=\varepsilon_{15}-\varepsilon_{13},
$$
and $\varepsilon_{45}=\varepsilon_{15}-\varepsilon_{14}$ (cf.~\eqref{eps'}).
We assume that the system is subject to the additional
moving affine constraint $\varepsilon_{15}=g(t)$. Combining all the constraints, we obtain
$$
\varepsilon= \varepsilon_{12} k_1 + \varepsilon_{13} k_2 + \varepsilon_{14} k_3 + g(t) k_0,
$$
where
$$
\begin{array}{ll}
k_1=(1,0,0,-1,-1,-1,0,0,0)'; & \ \ k_2=(0,1,0,1,0,0,-1,-1,0)';\\
k_3=(0,0,1,0,1,0,1,0,-1)'; & \ \ k_0=(0,0,0,0,0,1,0,1,1)'.
\end{array}
$$
This equation is equivalent to relation \eqref{eps} with
$
W={\rm span}\, (k_1,k_2,k_3).
$
The balance of forces \eqref{sig} at the nodes 2, 3, 4 reads
$$
\langle\sigma, k_1\rangle=\langle\sigma, k_2\rangle=\langle\sigma,k_3\rangle=0
$$
or, equivalently, $\sigma(t) \in W^\perp$ with $\sigma=(\sigma_{12},\sigma_{13},\sigma_{14},\sigma_{23},\sigma_{24},\sigma_{25},\sigma_{34},\sigma_{35},\sigma_{45})'$.
The matrix $A$ has the form
$$
A={\rm diag}\,(\alpha_{12},\alpha_{13},\alpha_{14},\alpha_{23},\alpha_{24},\alpha_{25},\alpha_{34},\alpha_{35},\alpha_{45}).
$$
In this example, the set $\Pi \cap V$  is is a 6-dimensional centrally symmetric convex polytope.

\subsection{The Prandtl-Ishlinskii model}\label{sect}
\begin{figure}[h]
	\begin{center}
		\includegraphics[width=4.3truecm]{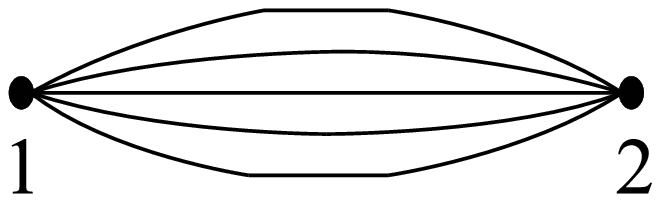} \quad \includegraphics[width=7truecm]{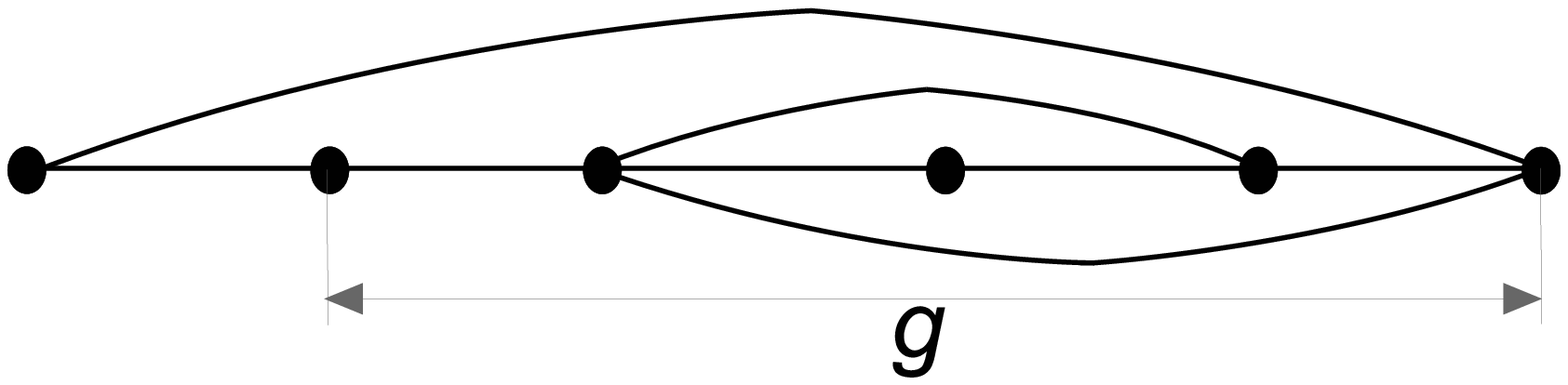}
	\end{center}
	\vskip-8truemm
	\centerline{(a) \hskip5.5truecm (b)\phantom{mmmmm}}
	%\vskip-2truemm
	\caption{\small
		(a) Parallel connection of springs in the Prandtl-Ishlinskii model.
		(b) A reducible connection of springs.
	}
	\label{fig4}
\end{figure}

The Prandtl-Ishlinskii operator ${\mathcal P}$ has scalar-valued inputs and outputs and is obtained as a weighted sum of a finite number of one-dimensional stop operators:
\begin{equation}\label{pi}
R={\mathcal P} [g]:=\sum_{n=1}^K \bar a_{n} \mathcal{S}_{\bar\rho_n}[g],
\end{equation}
where $g=g(t)$ and $R=R(t)$ are continuous input and output, respectively;
$\bar a_n$ is the Young's modulus of the $n$-th spring in the elastic domain; $\bar r_n:=\bar a_n\bar \rho_n$ is the maximal stress of the spring\footnote{The Prandtl-Ishlinskii model
	can include infinitely manu stops, in which case the sum in \eqref{pi} is replaced by the integral $R=\int_0^\infty a(\rho) \mathcal{S}_{\rho}[g] \,{\rm d\rho}$.
	However, in this work, only the finite (discretized) model \eqref{pi} is considered.}; and, without loss of generality we can assume that
\[
0<\bar \rho_1<\bar\rho_2<\cdots<\bar \rho_K.
\]
This operator is associated with the connection of Prandtl's springs shown in Fig.~\ref{fig4}(a), {i.e.} $\Gamma$ consists
of two vertices (nodes) $1$ and $2$ connected by
a multiedge. Formula \eqref{pi} relates the elongation
$
g=x_2-x_2-(x_2^*-x_1^*)
$
of the distance between the nodes (relative to the zero configuration) with the reaction
force $R$ applied at node $1$; the reaction at node $2$ equals $-R$.

\begin{figure}[h]
	\begin{center}
		\includegraphics[width=8truecm]{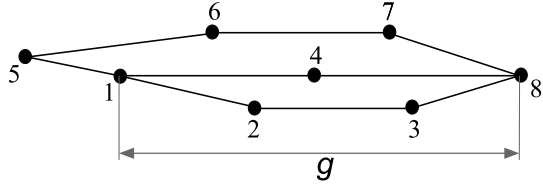}
	\end{center}
	\vskip-10truemm
	%\centerline{(a) \hskip5.5truecm (b)\phantom{mmmmm}}
	%\vskip-2truemm
	\caption{\small A ``linear'' connection of springs, which is equivalent to the Prandtl-Ishlinskii model.
	}
	\label{fig5}
\end{figure}

It is easy to see that the same relationship \eqref{pi} between the elongation $g$ of the moving constraint and the reaction force $R$ is valid for any simple connected graph $\Gamma$ of order $N$ %$N\ge 3$
in which every vertex $i=2,\ldots,N-1$ has degree $2$ and the moving affine constraint is applied
to the distance between the nodes $1$ and $N$ as in \eqref{input}.
Such connections will be called {\em linear}, see Fig.~\ref{fig5}.
%For such systems,
For linear connections, the parameters $K$, $\bar a_n$, $\bar \rho_n$  in \eqref{pi} are ``effective'' quantities that
can be related to the parameters $a_{ij}$, $\rho_{ij}$ of the springs (Young's moduli and yielding thresholds) via simple formulas.
Indeed, one can show that the stress $\sigma_{ij}=a_{ij} {\mathcal S}_{\rho_{ij}} [\varepsilon_{ij}]$ of each spring is related to the controlled distance $g=g(t)$ between the nodes $1$ and $N$
by another stop operator
$
\sigma_{ij} = {\tilde a_{ij}} {\mathcal S}_{\tilde \rho_{ij}} [g]
$
with appropriate effective $\tilde a_{ij}$ and $\tilde \rho_{ij}$.
For example, for the connection shown in Fig.~\ref{fig5},
$$
\tilde a_{12}=\tilde a_{23}=\tilde a_{38}=\frac{1}{\frac1{a_{12}}+\frac1{a_{23}}+\frac1{a_{38}}}, \qquad
\tilde a_{14}=\tilde a_{48}=\frac1{\frac1{a_{14}}+\frac1{a_{48}}},
$$
$$
\tilde a_{15}=\tilde a_{56}=\tilde a_{67}=\tilde a_{78}=\frac{1}{\frac1{a_{15}}+\frac1{a_{56}}+\frac1{a_{67}}+\frac1{a_{78}}}
$$
and
$$
\tilde r_{12}=\tilde r_{23}=\tilde r_{38}=\min\{r_{12},r_{23},r_{38}\},\quad
\tilde r_{14}=\tilde r_{48}=\min\{r_{14},r_{48}\},
$$
$$
\tilde r_{15}=\tilde r_{56}=\tilde r_{67}=\tilde r_{78}=\min\{r_{15},r_{56},r_{67},r_{78}\},
$$
where $r_{ij}=a_{ij}\rho_{ij}$, $\tilde r_{ij}=\tilde a_{ij}\tilde \rho_{ij}$.
Equations \eqref{sig''} now imply that $g$ is mapped to the reaction force $R$ of the constraint by the
Prandtl-Ishlinskii operator $R(t)= {\mathcal P} [g](t)=\tilde a_{38} {\mathcal S}_{\tilde\rho_{38}}[g]
+\tilde a_{48} {\mathcal S}_{\tilde\rho_{48}}[g]+\tilde a_{78} {\mathcal S}_{\tilde\rho_{78}}[g]$.
The general linear connection is considered in Section 4.

%\vskip-1truemm

An important characterization of Prandtl-Ishlinskii operator \eqref{pi} is the so-called {\em loading curve} $\phi: \mathbb R_+ \to \mathbb R_+$
defined as the response to the linear input:
\[
\phi={\mathcal P} [g_{id} ]\ \quad {\rm with} \ \quad g_{id}(t):=t, \ \ t\ge 0
\]
(where, in accordance with our agreement, we assume that $\sigma_n(0)=0$ for all the stop operators
$\sigma_n=\mathcal{S}_{\bar\rho_n}[g_{id}]$ in \eqref{pi}). Clearly, $\phi$ is a piecewise linear
function, which satisfies $\phi(0)=0$ and $\phi(\alpha)=const$ for sufficiently large $\alpha$. Equivalently,
\begin{equation}\label{vector}
\phi(0)=0,\qquad \frac{d\phi}{d\tau}(\tau)=\left\{
\begin{array}{rl}
\bar a_1+\bar a_2+\cdots+\bar a_K,& 0\le \tau<\bar\rho_1,\\
\bar a_2+\cdots+\bar a_K,& \bar\rho_1\le \tau<\bar\rho_2,\\
\vdots\phantom{m}\\
\bar a_K,& \bar\rho_{K-1}\le \tau<\bar\rho_K,\\
0, & \tau>\tau<\bar \rho_K.
\end{array}
\right.
\end{equation}
The loading curve uniquely identifies operator \eqref{pi} \cite{KP}, therefore we will write
${\mathcal P}={\mathcal P}_\phi$ when needed. Furthermore, for any $g=g(t)\in C$ with $g(0)=0$,
the output $R\in C$ defined by \eqref{pi} can be expressed explicitly in terms of the function $\phi$ as follows:
\begin{equation}\label{triangles}
R(t)=\phi(G_1(t)) + 2 \sum_{k\ge 2}  \phi\left(\frac{G_k(t)-G_{k-1}(t)}{2}\right),\quad \ t> 0,
\end{equation}
where the so-called {\em running main extremum values} $G_k$ of $g$
are defined by the relationships
\[
\begin{array}{l}
G_0(t):=\max \{|g(\tau)|:\ \tau\in[0, t]\}; \\ \tau_1:=\max\{\tau\in [0,t]:\ |g(\tau)|=G_0\};\quad G_1:=g(\tau_1)
\end{array}
\]
and, if $G_1\le 0$, then
\[
\begin{array}{l}
G_{2i}(t):=\max \{g(\tau):\ \tau\in[\tau_{2i-1},t]\}; \quad
\tau_{2i}:=\max\{\tau\in [0,t]:\ g(\tau)=G_{2i}\};\\
G_{2i+1}(t):=\min \{g(\tau):\ t\in[\tau_{2i}, t]\}; \quad
\tau_{2i+1}:=\max\{\tau\in[0,t]:\ g(\tau)=G_{2i+1}\};
\end{array}
\]
if $G_1> 0$, then
\[
\begin{array}{l}
\begin{array}{l}
G_{2i}(t):=\min \{g(\tau):\ \tau\in[\tau_{2i-1},t]\}; \quad
\tau_{2i}:=\max\{\tau\in [0,t]:\ g(\tau)=G_{2i}\};\\
G_{2i+1}(t):=\max \{g(\tau):\ t\in[\tau_{2i}, t]\}; \quad
\tau_{2i+1}:=\max\{\tau\in[0,t]:\ g(\tau)=G_{2i+1}\}
\end{array}
\end{array}
\]
for $i=1,2,\ldots$\
%the sequence of moments $\tau_k(t)$ satisfying $0<\tau_1(t)<\tau_2(t)<\cdots \le t$ at all times and
%$0<\tau_1(t)<\tau_2(t)<\cdots \le t$

By definition, a linear combination of operators \eqref{pi} with non-negative coefficients
is also a Preisach-Ishlinskii operator, $c_1{\mathcal P}_{\phi_1}+c_2{\mathcal P}_{\phi_2}={\mathcal P}_{c_1\phi_1+c_2\phi_2}$.
Remarkably, the class of the Prandtl-Ishlinskii operators \eqref{pi} is also closed with respect to composition due to the identity
${\mathcal P}_{\phi_1} \circ {\mathcal P}_{\phi_2}= {\mathcal P}_{\phi_1\circ \phi_2}$ \cite{bs,mroz}\footnote{The composition formula for Prandtl-Ishlinskii operators is based
	on Brokate's formula $({\rm Id}-{\mathcal S_{\rho_1}})\circ ({\rm Id}-{\mathcal S_{\rho_2}})=({\rm Id}-{\mathcal S_{\rho_1+\rho_2}})$ for the stop and play operators \cite{bs}.}.
Using these properties, one can associate a Prandtl-Ishlinskii operator
with more complex topologies $\Gamma$ than a linear connection such as in Fig.~\ref{fig5}.
In particular, for a multigraph with vertices $1,\ldots,N$, let us define an elementary graph operation that replaces a multiedge
with a simple edge; an elementary operation that replaces  a vertex $i\not\in\{1,N\}$ of degree $2$ and the two
edges emanating from it with
one edge; and, an elementary operation that eliminates a vertex $i\not\in\{1,N\}$ of degree $1$ and the edge emanating from it from the graph. We call a connected graph $G$ of order $N$ {\em reducible} if a finite number of such operations can reduce it
to the trivial graph that consists of two vertices $1$ and $N$ and an edge between them\footnote{Using analogy with electrical circuits
	of resistors, a reducible graph corresponds to a circuit that can be solved by applying parallel and series connection rules.}, see Fig.~\ref{fig4}(b).
One can show that if a connection of springs with a reducible graph $G$ is driven
by one moving constraint \eqref{input}, then the elongation $g$ of the constraint
is mapped to the stress of every spring by an (effective) Prandtl-Ishlinskii operator, $\sigma_{ij}={\mathcal P}_{\phi_{ij}}[g]$
and the reaction force is given by \eqref{pi} \cite{xx3}.

%Other simple connections of springs can be associated with the Prandtl-Ishlinskii operator. In particular, let us consider
%a connection of springs shown in Fig.~3.
%Hence, $g(t)$ is mapped to the total force $F(t)$ applied at node 1 by the.
%Prandtl-Ishlinskii operator $F(t)= {\mathcal P} [g](t)$.

\section{Main result}\label{main}
Consider a graph $\Gamma$ with edges $(ij)\in M$ and the corresponding system of $m=\# M$ connected Prandtl's springs with parameters $a_{ij}$ (stiffness) and $r_{ij}$ (maximal stress), see Sections \ref{1.1}, \ref{1.2}. Assume that the system is driven by one moving affine constraint
\eqref{input} and the set of all the constraints is described by inclusion  \eqref{eps}.

Denote by $u^*(t): [0, L]\to \Pi \cap V$ the solution of the differential inclusion \eqref{s}
with the input $g(t)=t$ and the initial value $u^*(0)=0$. The trajectory $\gamma$ of this solution
is a polyline $B_0 B_1\cdots B_\ell$ with one end $B_0$ at the origin. Each link of this polyline belongs to a different face
of the polytope $\Pi \cap V$, that is $B_{k-1} B_{k}\subset F_{k-1}$ for $k=1,\ldots,\ell$, where $F_k$ are (closed) faces of the polytope, which are all different, and $F_0=\Pi\cap V$.
Denote by $E_k$ the minimal affine subspace of $V$ which contains the face $F_k$ and by $\overset{\circ}{F_k}$ the interior of the set $F_k$ in $E_k$.

% denotes the interior of the polytope.
The parametrization $u^*(t)$ of the polyline $\gamma$ defines the partition
$0=d_0<d_1<\cdots <d_\ell=L$ of the interval $[0,L]$ by preimages of the points $B_k$: % links $B_{k-1} B_{k}$:
\begin{equation}\label{d}
u^*(0)=B_0=0,\quad u^*(d_1)=B_1,\quad u^*(d_2)=B_2, \quad \ldots, \quad u^*(d_\ell)=B_\ell.
\end{equation}

%Let us note that every face $F_k$ of the polytope $\Pi\cap V$ with $k\ge 2$ belongs to at least one $(m-1)$-dimensional face (facet)
%of the parallelepiped
%$$
%\Pi=\{y\in \mathbb{E}=\mathbb{R}^m:\ |y_{ij}|\le a_{ij}^{-\frac12} r_{ij} , \ (ij)\in M\}
%%%, \quad {\rm where} \quad \Delta_{ij}=a_{ij}^{-\frac12}=\rho_{ij} a_{ij}^{\frac12}.
%$$
%(cf.~\eqref{pi}). We introduce the notation
%$$
%\Pi_{ij}^+=\{y\in \Pi:\ y_{ij}=a_{ij}^{-\frac12} r_{ij} \},\qquad \Pi_{ij}^-=\{y\in \Pi:\ y_{ij}=-a_{ij}^{-\frac12} r_{ij}\}
%$$
%for the facets of $\Pi$.
%For every $(ij)\in M$, let us define a positive number $b_{ij}$ as follows:
%\begin{eqnarray}
%\label{b1}
%&& b_{ij}=L+\delta_{ij} \quad {\rm if} \quad B_1\not\in \Pi_{ij}^+ \cup \Pi_{ij}^-;\\
%\label{b2}
%&& b_{ij}=d_k \quad {\rm if} \quad B_k \in \Pi_{ij}^+ \cup \Pi_{ij}^- \quad {\rm and} \quad B_1,\ldots,B_{k-1} \not \in \Pi_{ij}^+ \cup \Pi_{ij}^-,
%\end{eqnarray}
%where $d_k$ is defined by \eqref{d} and $\delta_{ij}$ are arbitrary positive numbers.
\vskip2truemm

{\bf Theorem 1.} {\em
	Suppose that
	\begin{equation}\label{circ}
	F_{\ell-1} \subset F_{\ell-2}\subset\cdots\subset F_0; \qquad B_k\in \overset{\circ}{F_k} \ \ {\rm for} \ \ k=0,1,\ldots,\ell-1;
	\end{equation}
	and
	\begin{equation}\label{dim}
	{\rm dim}\, F_0 - {\rm dim}\, F_1 = {\rm dim}\, F_1 - {\rm dim}\, F_2 =\cdots ={\rm dim}\, F_{\ell-2} - {\rm dim}\, F_{\ell-1} = 1.
	\end{equation}
	Suppose that the parallelepiped
	\begin{equation}\label{omega}
	\Omega = \{ y = \tau_1 B_0B_1 +\cdots + \tau_\ell B_{\ell-1}B_\ell, \ |\tau_i|\le 1, \ i=1,\ldots,\ell\}
	\end{equation}
	belongs to the polytope $\Pi\cap V$.
	Assume that the map $u^*: [0,L]\to\gamma$ is invertible.
	Then, for every input $g(t)$ satisfying $g(0)=0$ and $|g(t)|\le L$ for all $t\ge 0$,
	the stress of each spring  relates to the variable $g$ via the %stop
	Preisach operator
	\begin{equation}\label{eff}
	%\sigma_{ij}(t) =\mu_{ij} \sqrt{a_{ij}}\, {\mathcal S}_{b_{ij}} [g](t),\qquad (ij)\in M,
	\sigma_{ij}(t) =\mathcal P_{\phi_{ij}}[g](t),\qquad (ij)\in M,
	\end{equation}
	with the loading curve
	%where $\mu_{ij}={\rm sign}\, f_0^{ij}$ with $f_0^{ij}$ denoting the $(ij)$-th coordinate of the vector $f_0$
	%defined by relationships \eqref{eps}, \eqref{f0};
	%and, the thresholds $b_{ij}$ are defined by \eqref{d}\,--%, \eqref{b1}, and \,\eqref{b2}.
	\begin{equation}\label{eff'}
	\phi_{ij}(\tau) =\sqrt{a_{ij}}\, u^*_{ij}(\tau), \qquad 0\le \tau\le L,
	\end{equation}
	where $u^*_{ij}$ is the $(ij)$-th component of the vector-valued function $u^*: [0,L]\to \Pi$.
}
\vskip2truemm

As we establish below, under the conditions of this theorem, the Moreau sweeping process \eqref{s}
behaves as an analog of the Prandtl-Ishlinskii operator with vector-valued outputs and the vector-valued $u^*: [0,L]\to \Pi\cap V$
acts as a counterpart of the loading curve for this operator. More precisely, the solution
of \eqref{s} with the zero initial condition $u(0)=0$ is given by the counterpart
of formula \eqref{triangles}:
%the output $R\in C$ defined by \eqref{pi} can be expressed explicitly in terms of the function $\phi$ as follows:
\begin{equation}\label{triangles'}
u(t)=u^*(G_1(t)) + 2 \sum_{i\ge 2}  u^*\left(\frac{G_i(t)-G_{i-1}(t)}{2}\right),\quad \ t> 0,
\end{equation}
where $G_k$ are running main extremum values of the scalar-valued input $g$; and, we extend the function $u^*$ to the interval
$[-L,L]$ by setting
\begin{equation}\label{odd}
u^*(-d)=-u^*(d),\qquad d\in [0,L].
\end{equation}
Formulas \eqref{eff}, \eqref{eff'} immediately follow from \eqref{triangles'} and the relationship
$\sigma=A^\frac12 u$.

\medskip
{\sl Proof of Theorem 1.} Given an input $g\in W^{1,1}([0,T];\mathbb{R})$ with $\max|g|\le L$,
we need to show that the function defined by \eqref{triangles'} satisfies $u(t)\in \Pi\cap V$ for all $t\in[0,T]$
and the inclusion \eqref{s} for almost every $t\in[0,T]$.
Suppose for definiteness that $G_1(t)\ge 0$.
Since
\begin{equation}\label{u*}
\dot u^*(d)=
\left\{
\begin{array}{cl}
B_0B_1/(d_1-d_0), & d_0< d<d_1,\\
B_1B_2/(d_2-d_1), & d_1<d < d_2,\\
\vdots & \\
B_{\ell-1}B_\ell/(d_\ell-d_{\ell-1}),& d_{\ell-1}<d<d_\ell,
\end{array}
\right.
\end{equation}
with $d_k$ defined in \eqref{d} (recall that $d_0=0$),
formula \eqref{triangles'} with $G_1(t)\ge 0$ implies
the relation
\begin{equation}\label{proof1}
u(t)=\sum_{k=1}^\ell \xi_k(t) \, \frac{B_{k-1}B_k}{d_k-d_{k-1}},\qquad t\ge 0,
\end{equation}
with
\begin{equation}\label{proof2}
\xi_k(t)=\bigl(\min \{G_1,d_k\}-d_{k-1}\bigr)^+ + \sum_{i\ge 2}(-1)^{i-1}\bigl(\min \{|G_i-G_{i-1}|, 2d_k\}- 2d_{k-1}\bigr)^+,
\end{equation}
where $a^+=\min\{a,0\}$ and $G_i=G_i(t)$. From the definition of the sequence $G_i$ given below Eq.~\eqref{triangles},
it follows that
\begin{equation}\label{compare}
2 G_1 \ge |G_2-G_1| \ge |G_3-G_2| \ge \cdots \ge |G_{i}-G_{i-1}| \ge \cdots
\end{equation}
where only a finite number of differences $G_{i}-G_{i-1}$ are nonzero. Therefore, \eqref{proof2} implies
\[
|\xi_k|\le d_k-d_{k-1},\qquad k=1,\ldots,\ell-1,
\]
and from \eqref{proof1} it follows that $u(t)$ belongs to the parallelepiped $\Omega$, which by assumption belongs to $\Pi\cap V$.
This proves $u(t)\in \Pi\cap V$ for the function \eqref{triangles'} in the case $G_1(t)\ge 0$. As the function $u^*$ in \eqref{triangles'} is odd,
it is easy to see that the same inclusion is valid for $G_1(t)< 0$.

It remains to prove \eqref{s}. As the solution operator of the Moreau sweeping process \eqref{g5}
subject to the zero initial condition $u(0)=0$ is continuous in the space $W^{1,1}([0,T];V)$
and so is the Prandtl-Ishlinskii operator in the space $W^{1,1}([0,T];\mathbb{R})$,
it suffices to show that function \eqref{triangles'} satisfies \eqref{s} for piecewise linear continuous inputs $g\in W^{1,1}([0,T];\mathbb{R})$.
Furthermore, given such an input, it suffices to establish that if the inclusion \eqref{s}
is valid for a.e.~$t$ from an interval $[0,\theta]\subset [0,T)$, then there is a $\delta>0$ such that this inclusion
is also true almost everywhere in $[\theta,\theta+\delta]$.

In order to establish \eqref{triangles'} is a solution of \eqref{s} on an initial small interval $[0,\theta]$, we
recall the {\em rate-independence} property
of the Moreau process and the Prandtl-Ishlinskii model,
which means that the input-output operator of each model commutes with increasing transformations of time \cite{bk}.
Let us take a sufficiently small $\theta>0$
so that $g$ is linear on $[0,\theta]$ and the trajectory $u^*(g(t))$, which satisfies $u^*(g(0))=u^*(0)=0$, remains in the interior of the polytope $\Pi\cap V$
for all $t\in [0,\theta]$. Then, \eqref{triangles'} is equivalent to $u(t)=u^*(g(t))$, which agrees with the rate-independence property of the Moreau process
for an increasing $g$. Since the function $u^*$ and the solution operator of the Moreau process are odd,
\eqref{triangles'} defines a solution of \eqref{s} on a sufficiently small interval $[0,\theta]$ for a decreasing $g$ too.

Assuming that \eqref{triangles'} satisfies \eqref{s} for $t\in[0,\theta]$, let $\hat G_k=G_k(\theta)$ be the sequence of the main extremum values of $g$ at the moment $\theta$.
%the sequence $\hat G_k$ is defined below Eq.~\eqref{triangles}.
Due to the central symmetry of $\Pi$, we can assume without loss generality that $\hat G_1>0$.
Since $g$ is piecewise linear, from the definition of $\hat G_k$ it follows that there is an integer $k_0$
such that $\hat G_k=g(\theta)$ for all $k\ge k_0$
and, depending on the parity of $k_0$, either
\begin{equation}\label{GG'}
\hat G_2 < \hat G_4< \cdots < \hat G_{2 p_0-2} < g(\theta)= \hat G_{2p_0-1} <\hat G_{2p_0-3}<\cdots< \hat G_3<\hat G_1
\end{equation}
with $k_0=2p_0-1$, or
\begin{equation}\label{GG}
\hat G_2 < \hat G_4< \cdots < \hat G_{2 p_0}=g(\theta) < \hat G_{2p_0-1} <\hat G_{2p_0-3} <\cdots< \hat G_3< \hat G_1
\end{equation}
with $k_0=2p_0$. The argument for both cases is similar, and we assume for definiteness that \eqref{GG'} holds.
%We consider these two cases separately.
%%%Note that if \eqref{GG} holds, then $g$ increases on $[t_j,t_{j+1}]$. In case \eqref{GG'}, $g$ decreases on $[t_j,t_{j+1}]$.

Since $g$ is piecewise linear, $\dot g=const$ on some interval $(\theta,\theta+\delta)$.
Using the rate-independence of the Moreau process, it suffices to consider the cases $\dot g=1$, $\dot g=-1$, and $\dot g=0$ for $t\in(\theta,\theta+\delta)$.
The latter case is trivial as $g=const$ implies that $u=const$ and $G_k=const$ in \eqref{triangles'} for all $k$ and all $t\in[\theta,\theta+\tau]$.
The cases when $g$ increases ($\dot g=1$) and $g$ decreases ($\dot g=-1$) on $(\theta,\theta+\tau)$ will be considered separately.

Let us introduce some notation. Suppose that the polytope $\Pi\cap V$ is the intersection of subspaces
\[
\bigcap_{i=1}^{m} \{ x\in V: \langle n_i, x\rangle \le c_i \} =\Pi\cap V,
\]
where each hyperplane $\langle n_i, x\rangle = c_i>0$ contains one facet of
$\Pi\cap V$ and $n_i$ is the unit outward normal vector to this facet.
Condition \eqref{dim} ensures that the face $F_k$ of $\Pi\cap V$
belongs to the intersection of $k$ hyperplanes $\langle n_i, x\rangle = c_i$,
hence we can number these hyperplanes in such a way that
\[
F_k \subset E_k=\bigcap_{i=1}^{k} \{ x\in V: \langle n_i, x\rangle = c_i \},\qquad k=1,\ldots,\ell-1.
\]
Recall that $\overset{\circ}{F}_k$ the interior of the face $F_k$ in the affine subspace
$E_k$ (with the agreement that $E_0=V$).
With this notation, the outward normal cone to $\Pi\cap V$ on $\overset{\circ}{F}_k$
%has the form
coincides with the positive linear span of the vectors $n_1,\ldots,n_k$:
\[
N_{\Pi\cap V}(x)=\left\{\sum_{i=1}^k \alpha_i n_i:\ \alpha_1,\ldots,\alpha_k\ge 0\right\},\qquad x\in \overset{\circ}{F}_k, \ \  k=1,\ldots,\ell-1.
%\{y\in V: \langle y, n_i\rangle\ge 0,\, i=1,\ldots,k\},\qquad x\in \overset{\circ}{F}_k, \ \  k=1,\ldots,\ell-1.
\]
%Furthermore,
%denoting by $E_k^o$ the orthogonal complement to the span of vectors $n_1,\ldots, n_{k-1}$ in $V$:
%\[
%E_k^o:=\bigcap_{i=1}^{k} \{ x\in V: \langle n_i, x\rangle = 0 \},\qquad k=1,\ldots,\ell-1,
%\]
%and by $f_k$ the orthogonal projection of the vector $f_0$ onto $E_k^o$:
%\[
%f_k={\rm proj}_{E_k^o}\,f_0
%\]
%(where, again, $E_0^o=V$), one can see that the inclusion \eqref{s} for $u\in \overset{\circ}{F}_k$ is equivalent to the equality
%\[
%\dot u = \dot g f_k, \qquad u\in \overset{\circ}{F}_k,\quad k=0,\ldots,\ell-1.
%\]
%In particular, links of the polyline $\gamma$ satisfy
%\[
%B_{k-1}B_k = (d_k-d_{k-1}) f_{k-1},\qquad k=1,\ldots,\ell.
%\]
In particular, since relations \eqref{circ} ensure that $u^*(d)\in \overset{\circ}{F}_{k-1}$ for $d_{k-1}<d<d_k$,
the inclusion \eqref{s} for $\dot u^*$ is equivalent to
\[
-\dot u^*(d) + f_0\in \left\{\sum_{i=1}^{k-1} \alpha_i n_i:\ \alpha_1,\ldots,\alpha_{k-1}\ge 0\right\} \quad {\rm for} \quad d_{k-1}<d<d_k
\]
with $k=1,\ldots,\ell-1$. Combining these relations with \eqref{u*}, we obtain
\begin{equation}\label{fin}
-\frac{ B_{k-1} B_{k}}{d_k-d_{k-1}} + f_0 \in \left\{\sum_{i=1}^{k-1} \alpha_i n_i:\ \alpha_1,\ldots,\alpha_{k-1}\ge 0\right\}.
\end{equation}
One can also see that the inclusion \eqref{s} for $u\in \overset{\circ}{F}_k$
is equivalent to the equality
\begin{equation}\label{cir}
\dot u = \dot g f_k,
\end{equation}
where
$$
f_k={\rm proj}_{\overset{\circ}{E}_k} f_0
$$ is the orthogonal projection of the vector $f_0$ onto the
%orthogonal complement $E^\perp_k=
subspace
$\overset{\circ}{E}_k=
\{ x\in V: \langle n_i, x\rangle = 0, \ i=1,\ldots,k\}$, which is
%of the affine subspace $E_k$ in $V$.
parallel to $E_k$ in $V$. Note that since $U=V^\perp$,
$$
f_k={\rm proj}_{\overset{\circ}{E}_k}A^\frac12 k_0.
$$

Now, suppose that \eqref{GG'} holds at the moment $\theta$ and $\dot g=-1$ for $t\in(\theta,\theta+\delta)$.
If $\delta>0$ is sufficiently small, then according to the definition of the main extremum values of $g$,
one has
\[
G_2(t) < \cdots < G_{2p_0-2}(t) < G_{2 p_0}(t)=g(t) < G_{2p_0-1}(t) <G_{2p_0-3}(t) <\cdots< G_1(t)
\]
with
\[
G_1(t)=G_1(\theta),\ G_2(t)=G_2(\theta),\ \ldots,\ G_{2p_0-1}(t)=G_{2p_0-1}(\theta) \ \ {\rm for\ all} \ \ t\in[\theta,\theta+\delta]
\]
(cf.~\eqref{GG}). Therefore, comparing formula \eqref{triangles'} at the moments $\theta$ and $t\in (\theta,\theta+\delta)$ (note that the sum in
\eqref{triangles'} contains nonzero terms for $2\le k\le 2p_0-1$ at the moment $\theta$ and
for $2\le k\le 2p_0$ at the moment $t$), we obtain
\[
u(t)=u(\theta)+2 u^*\left(-\frac{(t-\theta)f_0 }2\right),
\]
and using \eqref{odd}, \eqref{u*} and $\delta<d_1$, we arrive at the relation
\[
\dot u = -f_0,\qquad t\in[\theta,\theta+\delta],
\]
which agrees with \eqref{s}.

Finally, suppose that $\dot g=1$ for $t\in(\theta,\theta+\delta)$ and, again,
relations \eqref{GG'} are valid at the moment $\theta$.
In this case, assuming that $\delta>0$ is sufficiently small,
\[
G_2(t) < \cdots < G_{2p_0-2}(t) < g(t) = G_{2p_0-1}(t) <G_{2p_0-3}(t) <\cdots< G_1(t)
\]
with
\[
G_1(t)=G_1(\theta),\ G_2(t)=G_2(\theta),\ \ldots,\ G_{2p_0-2}(t)=G_{2p_0-2}(\theta) \ \ {\rm for\ all} \ \ t\in[\theta,\theta+\delta].
\]
Therefore, \eqref{triangles'} implies
\[
u(t)=u^*(\theta)-2u^*\left(\frac{G_{2p_0-1}(\theta)-G_{2p_0-2}(\theta)}{2}\right)+2u^*\left(\frac{g(t)-G_{2p_0-2}(\theta)}{2}\right),
\]
hence
\begin{equation}\label{dot}
\dot u(t)=\dot u^*\left(\frac{g(t)-G_{2p_0-2}(\theta)}{2}\right),\qquad t\in(\theta,\theta+\delta).
\end{equation}
Now, note that there is a $k$, $1\le k\le\ell$, such that
\[
2d_{k-1}\le G_{2p_0-1}(\theta)-G_{2p_0-2}(\theta) <2d_k,
\]
and since $g(t)>G_{2p_0-2}(\theta)$, one also has
\begin{equation}\label{xx}
2d_{k-1}< g(t)-G_{2p_0-2}(\theta) <2d_k,\qquad t\in(\theta,\theta+\delta),
\end{equation}
for a sufficiently small $\delta$. Combining these relations with \eqref{u*} and \eqref{dot}, we see that
\begin{equation}\label{dot'}
\dot u(t)=\frac{B_{k-1}B_{k}}{d_k-d_{k-1}},\qquad t\in(\theta,\theta+\delta).
\end{equation}
Also, relations \eqref{compare} and \eqref{xx} imply
\[
2 G_1\ge |G_2-G_1| \ge \cdots \ge |g(t)-G_{2p_0-2}| > 2d_{k-1}\ge 2d_j,\qquad j=1,\ldots, k-1,
\]
where $G_i=G_i(\theta)=G_i(t)$ for $i\le 2p_0-2$ and $g(t)=G_{2p_0-1}(t)$.
Hence, for $j=1,\ldots, k-1$, \eqref{proof2} becomes
\[
\xi_j(t)=d_j-d_{j-1} +2\sum_{i=2}^{2p_0-1} (-1)^{i-1}(d_j-d_{j-1})=d_j-d_{j-1}, \qquad t\in(\theta,\theta+\delta),
\]
(note that all terms with $i\ge 2p_0$ in the sum in \eqref{proof2} are zero) and \eqref{proof1} implies
\[
u(t)=B_{k-1} + \sum_{j=k}^\ell \xi_k(t)\frac{B_{k-1}B_k}{d_k-d_{k-1}}, \qquad t\in(\theta,\theta+\delta).
\]
Since $B_{k-1},\ldots,B_\ell \in F_{k-1}$, it follows that $u(t)\in F_{k-1}$ and hence
\[
n_1,\ldots,n_{k-1}\in N_{\Pi\cap V}(u(t)), \qquad t\in(\theta,\theta+\delta).
\]
Combining this relation with \eqref{fin} and \eqref{dot'}, we conclude that \eqref{s} with $\dot g=1$ holds for
$t\in(\theta,\theta+\delta) $. This completes the proof.

%The other case when and $\dot g=1$ for $t\in(\theta,\theta+\delta)$
%leads to two possibilities depending on ...
%We assume again that \eqref{GG'} is valid at the moment $\theta$.

\section{Discussion}\label{discussion}
A few remarks are in order.

{\bf  1.}
Under the assumptions Theorem 1, the complex connection of springs defined by the graph $\Gamma$ responds
to arbitrary variations $g(t)$ of the given moving affine constraint in the same way as a simple parallel connection
of springs described by an effective Prandtl-Ishlinskii opearator. Indeed, the stress of each spring
of the complex system is related to $g(t)$ via the Prandtl-Ishlinskii operator \eqref{eff}
and the total force applied to the system at node 1 relates to the displacement $g$  via the Prandtl-Ishlinskii operator (cf.~\eqref{sig''})
$$
-R=\sum_{(1j)\in M} \sigma_{1j} = \sum_{(1j)\in M} \mathcal P_{\phi_{1j}}[g](t).
$$
\vskip2truemm

{\bf  2.} %Fig.~?? illustrates the geometry of hysteresis loops produced by the Moreau model under the conditions of Theorem 1.
According to formula \eqref{triangles'}, %these loops
under the conditions of Theorem 1,
hysteresis loops of the Moreau model and the loops of the Prandtl-Ishlinskii model
have similar properties and can be constructed by the same simple manipulations with the graph of the loading curve, see \cite{xx3} for details. The only difference is that this graph
is two-dimensional for the Prandtl-Ishlinskii model, while the loading curve $\gamma$ for the Moreau model in Theorem 1 is multi-dimensional.

\begin{figure}[h]
	\begin{center}
		\includegraphics[width=5.5truecm]{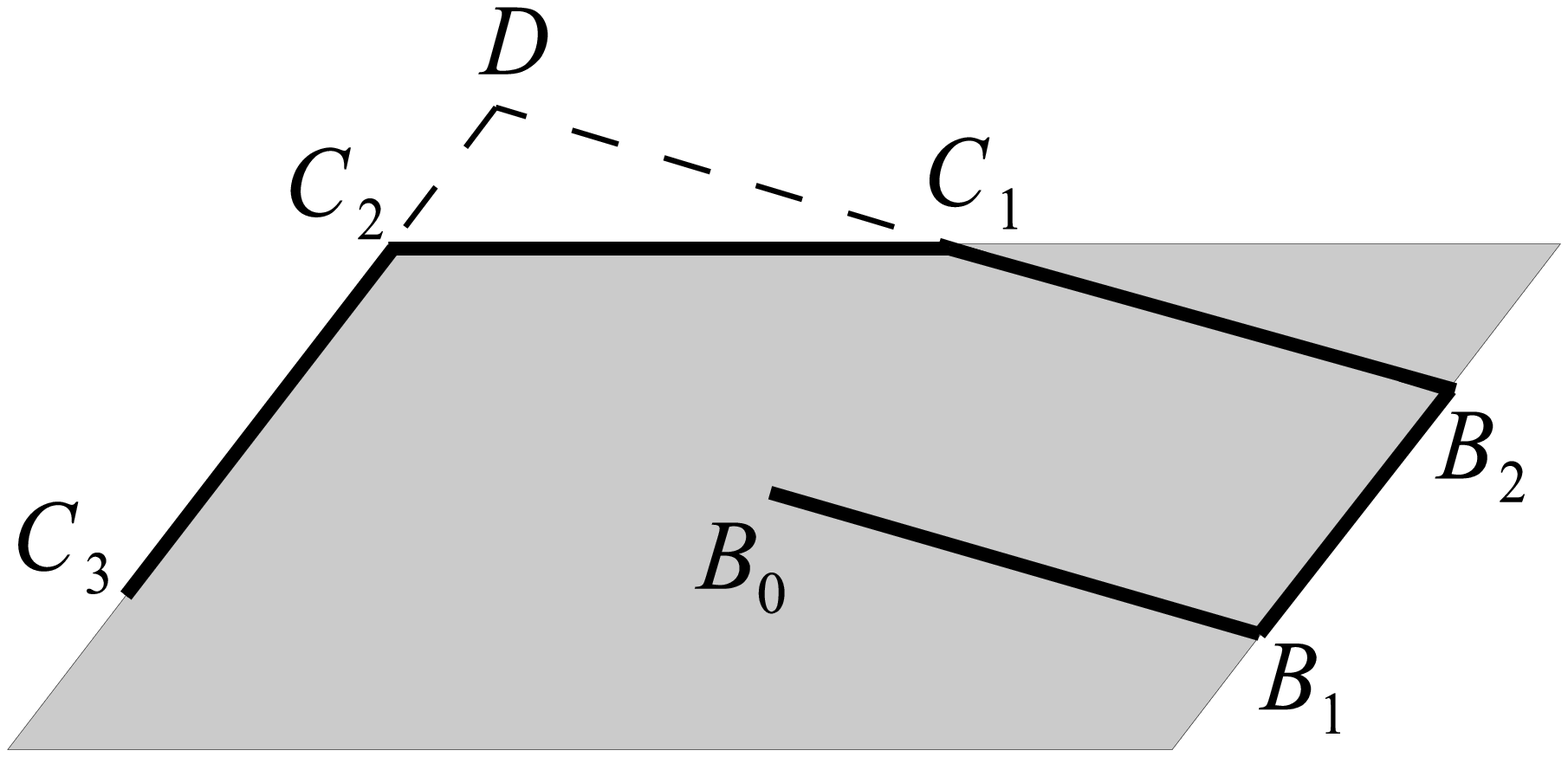} \ \ \includegraphics[width=5.5truecm]{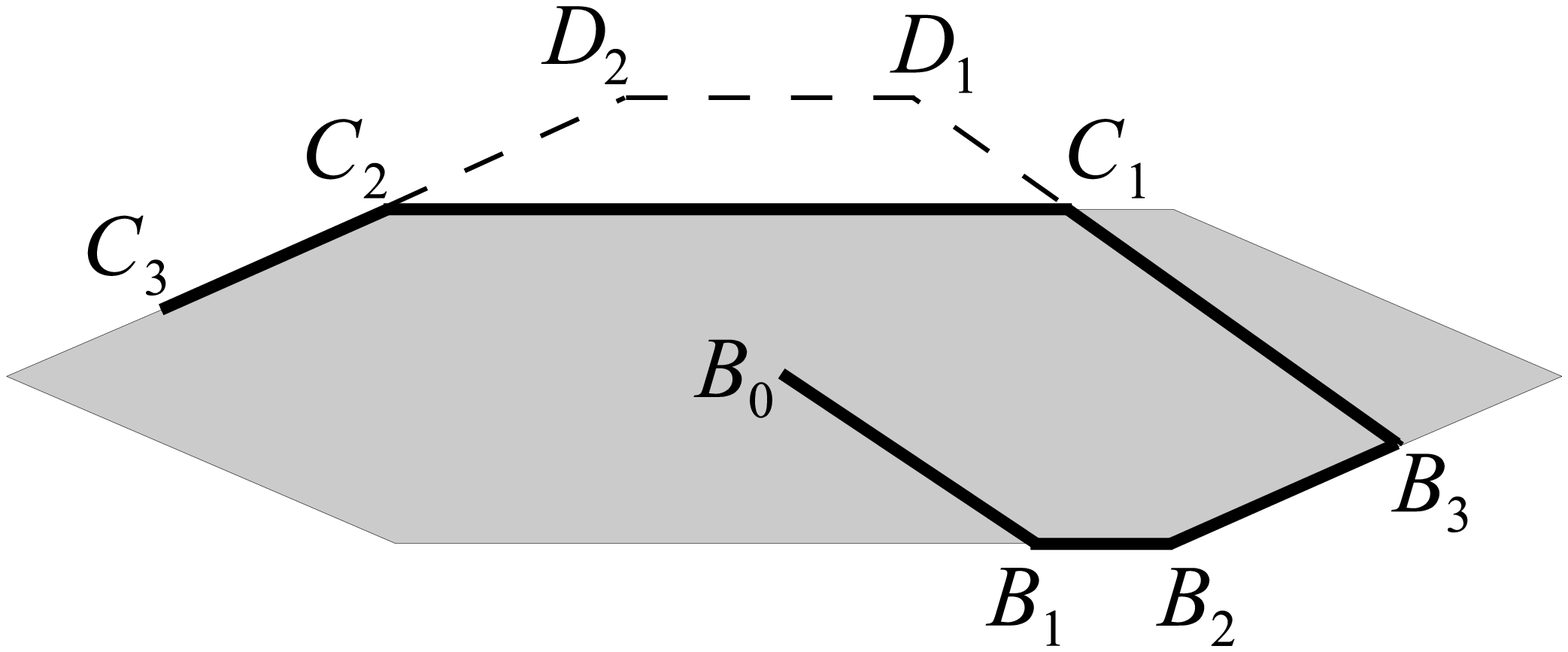}
	\end{center}
	\vskip-8truemm
	\centerline{\hskip1truecm (a) \hskip5.5truecm (b)\phantom{mmmmm}}
	%\vskip-2truemm
	\caption{\small Violations of conditions of Theorem 1.
		The shaded area represents the polytope $\Pi\cap V$.
		The vector $f_0$ used in the definition \eqref{s} of the Moreau process points in the direction of the vector $B_0B_1$.
		(a) The prism $\Omega$ does not belong to $\Pi\cap V$. The polyline $B_0B_1B_2C_1C_2C_3$
		that contains the the polyline $\gamma=B_0B_1B_2$
		represents the trajectory of a solution $u$ to \eqref{s} for the input $g(t)f_0$
		where $g$ increases from zero to a maximum value $g_*$ and then decreases to the minimum value $-g_*$.
		Formula \eqref{triangles'} defines a different polyline $B_0B_1B_2DC_3$ for the same input.
		(b) Conditions \eqref{circ}--\eqref{dim} are violated because ${\rm dim}\, F_{1} = {\rm dim}\, F_2$.
		The polylines $B_0B_1B_2B_3C_1C_2C_3$ and $B_0B_1B_2B_3D_1D_2C_3$ the trajectory of
		the inclusion \eqref{s} and curve prescribed by formula \eqref{triangles'} in response to the same input as in panel (a).
	}
	\label{fig6}
\end{figure}

\vskip2truemm
{\bf 3.} Fig.~\ref{fig6} illustrates the role of conditions $\Omega\subset \Pi\cap V$ and ${\rm dim}\, F_{k-1} - {\rm dim}\, F_k=1$ of Theorem 1. Fig.~\ref{fig6}(a) violates condition \eqref{omega} because
$\Omega\not\subset \Pi\cap V$.
In Fig.~6(b), ${\rm dim}\, F_{1} = {\rm dim}\, F_2$, hence conditions
\eqref{circ}--\eqref{dim} are violated. In both cases, we see that the trajectory of the Moreau process
does not have the shape prescribed by formula \eqref{triangles'}, which generalizes the Prandtl-Ishlinskii operator.
%are not of Prandtl-Ishlinskii type, because they lose the symmetry property described in Remark 1.

%\vskip-5truemm

\vskip2truemm
{\bf 4.} The assumption that the map $u^*: [0,L]\to\gamma$ is invertible is made for simplicity.
It is straightforward to extend the theorem to the case when it is not satisfied.

\begin{figure}[h]
	\begin{center}
		\includegraphics[width=4.5truecm]{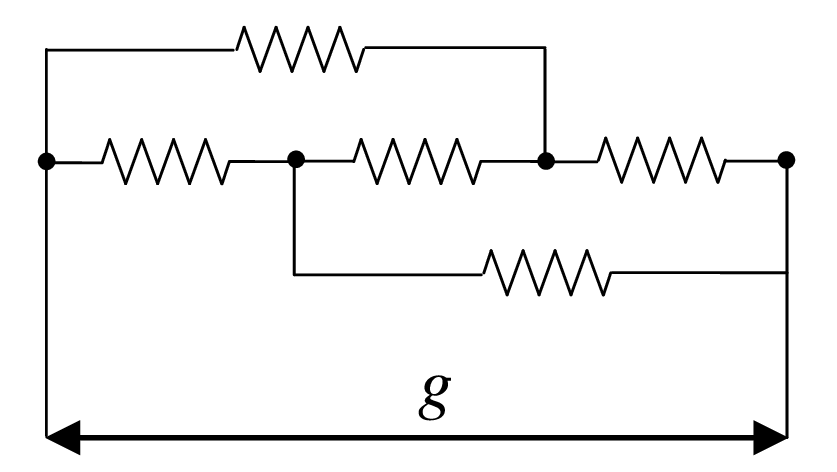} 
	\end{center}
	\caption{A `minimal' system, which may be not reducible to the PI operator (depending on the parameters of springs).}
	\label{5springs}
\end{figure}

\vskip2truemm
{\bf 5.}
Examples in which the system of springs shown in Fig.~\ref{fig1} cannot be reduced to a Prandtl-Ishlinskii operator
can be found in \cite{xx3}. %{\bf Add an explicit example?}
The simplest example is presented in Fig.~\ref{5springs}. It is easy to show that this system of 5 springs is or is not be reducible to the PI operator (according to Theorem 1) depending on the parameters of the springs.

%\vskip-8truemm
\begin{figure}[h]
	\begin{center}
		\includegraphics[width=8truecm]{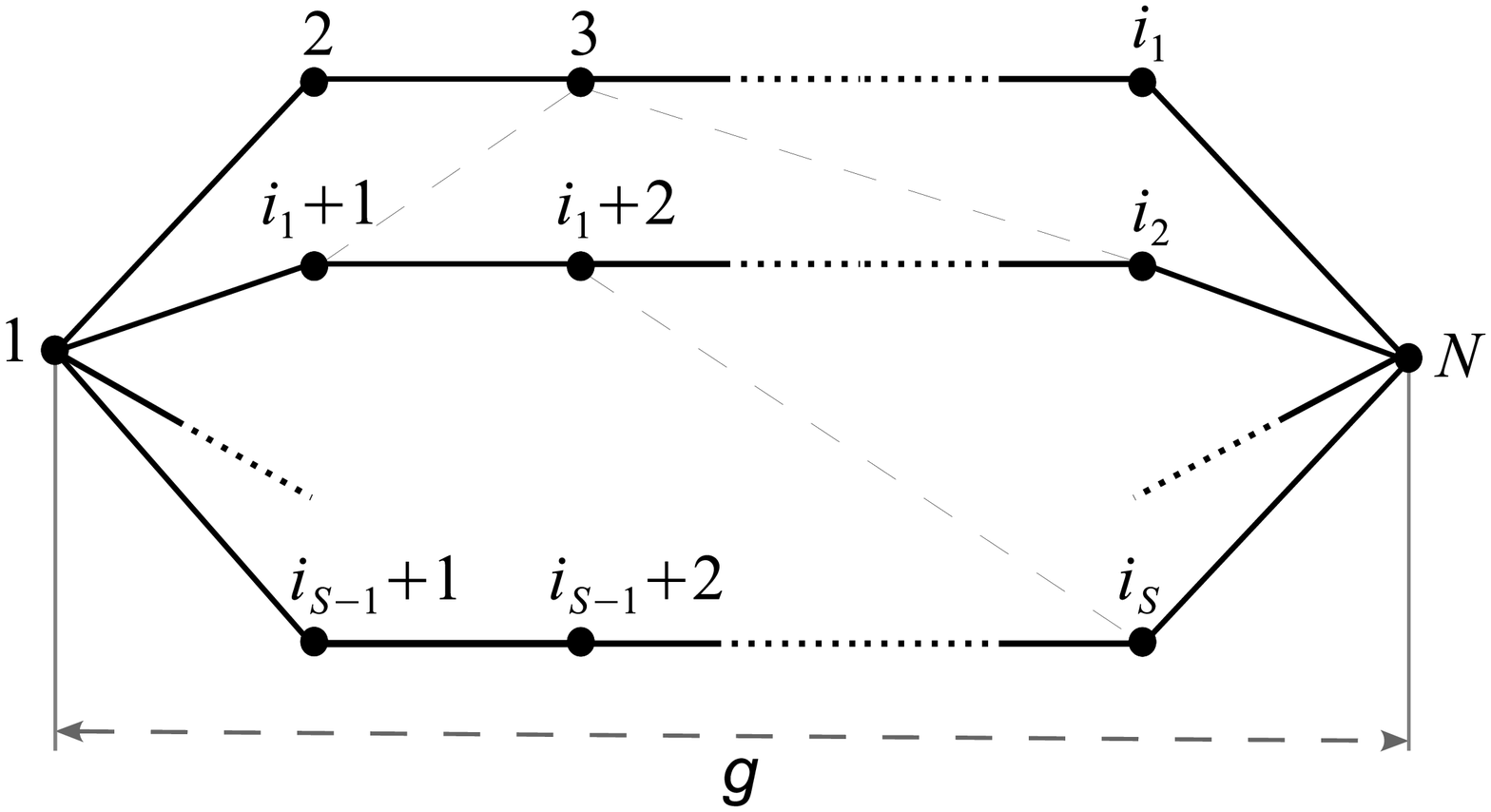}
	\end{center}
	\vskip-11truemm
	\caption{\small
		A graph $G$ of a general {\em linear} connection of $m$ springs (thick lines).
		Thin dashed lines represent possible additional edges as introduced in Theorem 2.
	}
	\label{fig7}
\end{figure}

%\vskip2truemm
%\vskip2truemm
%{\bf Remark 5.}
%
%\vskip2truemm
{\bf 6.}
Examples with complex topology that satisfy the conditions of Theorem 1 can be created by a perturbation of simple systems considered in Section \ref{sect}.
In particular, let us consider a {\em linear} connection of springs with the corresponding graph $\Gamma$ shown in Fig.~\ref{fig7} (the graph
shown in Fig.~\ref{fig5} is an example of such a graph). Here
\[
1=i_0<i_1<\cdots<i_{S-1}<i_S=N-1
\]
and the list of all edges $(ij)$ with $i<j$ is as follows:
node $1$ is connected with nodes $\{i_0+1,i_1+1,\ldots,i_{S-1}+1\}$;
node $N$ is connected with nodes $\{i_1,i_2,\ldots,i_S\}$;
and, each node $i$ with $i_n+1\le i<i_{n+1}$ is connected with node $i+1$
for $0\le n\le S-1$. The total number of edges equals $m=i_S+S-1$.
Denote
\begin{equation}\label{ta}
\tilde a_n:=\left({\frac1{a_{1,i_n+1}}+\sum_{i=i_n+1}^{i_{n+1}-1}\frac1{a_{i,i+1}}+
	\frac1{a_{{i_{n+1},N}}}}\right)^{-1},%\qquad n=0,\ldots,S-1.
\end{equation}
\begin{equation}\label{tr}
\tilde r_n:=\min \{r_{1,i_n+1},r_{i_n+1,i_n+2},\ldots,r_{i_{n+1}-1,i_{n+1}},r_{i_{n+1},N}\}
\end{equation}
with $n=0,\ldots,S-1$ (see example in Section 2.6) and set
\[
q_n(i)=
\left\{
\begin{array}{lll}
i_n+1, & i=1,\\
i+1, & i_n\le i\le i_{n+1}-1,\\
N, & i_{n+1}.
\end{array}
\right.
\]

\vskip2truemm
{\bf Theorem 2.} {\em
	Consider a linear connection of $m$ springs (see Fig.~\ref{fig7}). Assume that
	\begin{equation}\label{f}
	r_{i,q_n(i)} \ne r_{j,q_n(j)},\qquad i\ne j; \ \ i,j\in \{1, i_n+1,i_n+2,\ldots, i_{n+1}\},
	\end{equation}
	\begin{equation}\label{ff}
	\frac{\tilde r_{n_1}}{\tilde a_{n_1}}\ne \frac{\tilde r_{n_2}}{\tilde a_{n_2}} ,\qquad n_1\ne n_2
	\end{equation}
	for $0\le n, n_1,n_2\le S-1$.
	Let us extend this system with any set of additional connections of nodes %by {\em ideal}
	by springs with thresholds $r_{ij}$ and Young's moduli $a_{ij}$.
	Then, given any $r_{th}>0$ there is a $\delta=\delta(r_{th})>0$ such that
	if for all the added connections
	$r_{ij}\ge r_{th}$ and simultaneously $a_{ij}<\delta$, i.e.~the Young's moduli $a_{ij}$ of the added springs are sufficiently small,
	then the extended system of springs satisfies all the conditions of Theorem 1.
}

\vskip2truemm
%{\bf Remark 5.} An ideal spring can be viewed as Prandtl's spring with an infinitely large maximal value $\rho_{ij}=r_{ij}/a_{ij}=\infty$ of elastic deformation.
%Naturally, given any $r_0>0$, the ideal springs in the formulation of Theorem 2 can be replaced with Prandtl's springs
%that all have yielding thresholds satisfying $r_{ij}\ge r_0>0$ and sufficiently small Young's moduli $a_{ij}$.
%The conclusion of the theorem remains valid.
%{\bf $S$ should be $\ell-1$?}

\vskip2truemm
{\sl Proof of Theorem 2.} First, we consider the linear connection of $m$ springs.
For this connection, %the subspace $W\subset \mathbb{R}^m$ and the vector $k_0$ in
the set of geometric and moving affine constraints \eqref{e1} has the form
\[
\varepsilon_{1,i_n+1}+\sum_{i=i_n+1}^{i_{n+1}-1}\varepsilon_{i,i+1}+\varepsilon_{i_{n+1},N}=g,\qquad\quad n=0,\ldots,S-1.
\]
Therefore, the subspace $W\ni \varepsilon=e+p$ of $\mathbb{E}=\mathbb{R}^m$ is defined by the relations
\[
\varepsilon_{1,i_n+1}+\sum_{i=i_n+1}^{i_{n+1}-1}\varepsilon_{i,i+1} +\varepsilon_{i_{n+1},N}=0,\qquad\quad n=0,\ldots,S-1,
\]
its orthogonal complement $W^\perp\ni \sigma$ is given by
\[
\sigma_{1,i_n+1}=\sigma_{i_n+1,i_n+2}=%\sigma_{i_n+2,i_n+3}=
\cdots=\sigma_{i_{n+1}-1,i_{n+1}}=\sigma_{i_{n+1},N},\qquad\quad n=0,\ldots,S-1,
\]
and the components $k_0^{i,j}$ of the vector $k_0$ in \eqref{e1} have the form
\[
k_0^{i,j}=\left\{ \begin{array}{ll}
1, & i=i_n, \ j=N,
\\
0, & {\rm otherwise},
\end{array}\right.\qquad n=1,\ldots,S.
\]
Therefore, the orthogonal subspaces $U=A^\frac12 W\ni u$ and $V=A^{-\frac12}W^\perp=U^\perp\ni v$ (see \eqref{UV}) are defined by
the systems
\[
\frac{u_{1,i_n+1}}{\sqrt{a_{1,i_n+1}}}+
\sum_{i=i_n+1}^{i_{n+1}-1}
\frac{u_{i,i+1}}{\sqrt{a_{i,i+1}}}+
\frac{u_{i_{n+1},N}}{\sqrt{a_{{i_{n+1},N}}}}=0,
\]
\begin{equation}\label{vvvv}
v_{1,i_n+1}{\sqrt{a_{1,i_n+1}}}=v_{i_n+1,i_n+2}{\sqrt{a_{i_n+1,i_n+2}}}=
%=v_{k_n+2,k_n+3}
%=\cdots=v_{k_{n+1}-1,k_{n+1}}{\sqrt{a_{k_{n+1}-1,k_{n+1}}}}
\cdots=v_{i_{n+1},N}{\sqrt{a_{{i_{n+1},N}}}}
\end{equation}
with $n=0,\ldots,S-1$, respectively.
From these relations, it follows that the components of vector \eqref{f0} are given by
\begin{equation}\label{add1}
f_0^{1,i_n+1}{\sqrt{a_{1,i_n+1}}}=f_0^{i_n+1,i_n+2}{\sqrt{a_{i_n+1,i_n+2}}}=
%=v_{k_n+2,k_n+3}
%=\cdots=v_{k_{n+1}-1,k_{n+1}}{\sqrt{a_{k_{n+1}-1,k_{n+1}}}}
\cdots=f_0^{i_{n+1},N}{\sqrt{a_{{i_{n+1},N}}}}=\tilde a_n,
\end{equation}
and the polytope $\Pi\cap V$ is an $S$-dimensional parallelepiped
defined by relations \eqref{vvvv} and
\begin{equation}\label{fff}
|v_{1,i_n+1}|{\sqrt{a_{1,i_n+1}}}
%=|v_{i_n+1,i_n+2}|{\sqrt{a_{i_n+1,i_n+2}}}= \cdots=|v_{i_{n+1},N}|{\sqrt{a_{{i_{n+1},N}}}}
\le \tilde r_n,
\end{equation}
where $\tilde a_n$, $\tilde r_n$ are defined in \eqref{ta}, \eqref{tr}.
Using \eqref{f}, \eqref{ff}, without loss of generality, we can assume that
\begin{equation}\label{z}
\tilde \rho_0<\tilde\rho_1<\cdots<\tilde \rho_{S-1},\qquad \tilde\rho_n:=\frac{\tilde r_n}{\tilde a_n},
\end{equation}
and $\tilde r_n=r_{1,i_n+1}$ for each $n=0,\ldots,S-1$, i.e.
\begin{equation}\label{zz}
\tilde r_n=r_{1,i_n+1}< \min\{r_{i_n+1,i_n+2},\ldots,r_{i_{n+1}-1,i_{n+1}},r_{i_{n+1},N}\},
\end{equation}
(because the edges can always be labeled so as to ensure these relationships).
Now, from \eqref{cir}, \eqref{vvvv} and \eqref{fff}, it is easy to derive explicit expressions for the components of vectors $f_1,\ldots,f_{S}$
and coordinates of points $B_1,\ldots,B_S$:%,B_{S+1}$:
\footnote{One can also use explicit balance of forces to obtain the same result.}
\begin{equation}\label{add2}
f_j^{i,q_n(i)}=
\left\{
\begin{array}{cl}
0, & i\in \{1, i_n+1,i_n+2,\ldots, i_{n+1}\}, \ n< j,\\
f_0^{i,q_n(i)}, & i\in \{1, i_n+1,i_n+2,\ldots, i_{n+1}\}, \ n\ge j,
\end{array}
\right.
\end{equation}
\[
B_j^{1,i_n+1}{\sqrt{a_{1,i_n+1}}}=B_j^{i_n+1,i_n+2}{\sqrt{a_{i_n+1,i_n+2}}}=
%=v_{k_n+2,k_n+3}
%=\cdots=v_{k_{n+1}-1,k_{n+1}}{\sqrt{a_{k_{n+1}-1,k_{n+1}}}}
\cdots=B_j^{i_{n+1},N}{\sqrt{a_{{i_{n+1},N}}}}=\tilde a_n \min\{\tilde \rho_n,\tilde \rho_{j-1}\}
\]
with $n=0,\ldots,S-1$.

These explicit formulas ensure that:

\begin{itemize}
	
	\item[(i1)]
	The point $B_j\in V$ belongs to the open face $\overset{\circ}{(\Pi_1\cap \cdots \cap \Pi_j)}$ of the parallelepiped $\Pi$, i.e.
	\[
	B_j\in \overset{\circ}{(\Pi_1\cap \cdots \cap \Pi_j)},
	\]
	where $\Pi_j$ denotes the facet $$\Pi_j=\{v\in \Pi: v_{1,i_{j-1}+1}=r_{{1,i_{j-1}+1}}/\sqrt{a_{1,i_{j-1}+1}}\}$$ of $\Pi$ for $j=1,\ldots,S$.
	
	\item[(i2)]
	$V\cap \Pi_1\cap \cdots \cap \Pi_j =F_j$ and conditions \eqref{circ}, \eqref{dim} are satisfied.
	
	\item[(i3)]
	The subspace $V$ and the hyperplanes
	$$
	\Xi_j=\{v\in \mathbb{E}: v_{i_{j-1}+1}=0\}, \qquad j=1,\ldots,S,
	$$
	which are parallel to the facets $\Pi_j$ of $\Pi$, are in general linear position in the sense that
	\[
	{\rm dim} (V\cap \Xi_1\cap\cdots\cap \Xi_j)={\rm dim (V)}-j=S-j,\qquad j=1,\ldots,S.
	\]
	This implies $\overset{\circ}{E}_j=V\cap \Xi_1\cap\cdots\cap \Xi_j$ in the notation of the previous section.
	
	\item[(i4)] $f_j={\rm proj}_{V\cap\Xi_0\cap \Xi_1\cap\cdots\cap \Xi_j} A^\frac12 k_0\ne 0$ for $j=0,\ldots,S-1$ (where $\Xi_0=\mathbb{E}$). 
	
\end{itemize}

The next step is to show that the same properties (i)-(iv) are satisfied for the extended system of $\tilde m$ springs.
To this end, we first consider a simple embedding of the geometric picture considered above into the space $\tilde{\mathbb{E}}=\mathbb{R}^{\tilde m}$
by identifying each vector $v\in \mathbb{E}=\mathbb{R}^m$ with components $v_{ij}$, $(ij)\in M$, with the vector $\tilde v\in \tilde{\mathbb{E}}$
with the components
\[
\tilde v_{ij}=
\left\{
\begin{array}{cl}
v_{ij}, & (ij)\in M,\\
0, & (ij)\in \tilde M\setminus M.
\end{array}
\right.
\]
(recall that $\Gamma$ is a subgraph of $\tilde \Gamma$, hence $M\subset \tilde M$). With this identification,
$\mathbb{E}$ becomes a subspace $\tilde{\mathbb{E}}$ and we denote by $\mathbb{E}^\perp$ its orthogonal complement
with respect to the scalar product
\[
\langle y,z\rangle =\sum_{(ij)\in \tilde M} y_{ij}z_{ij},\qquad y,z\in\tilde{\mathbb{E}},
\]
(cf.~\eqref{sp}), i.e. $\mathbb{E}^\perp=\{ v\in \tilde{\mathbb{E}}: v_{ij}=0, (ij)\in \tilde M\setminus M\}$. Also, define the sets
\[
\tilde V=V\oplus \mathbb{E}^\perp,\quad \tilde \Pi=\Pi\oplus \mathbb{E}^\perp, \quad \tilde \Pi_j=\Pi_j\oplus \mathbb{E}^\perp,\quad
\tilde \Xi_j=\Xi_j\oplus \mathbb{E}^\perp
\]
and the diagonal $\tilde m\times \tilde m$ matrix $\tilde A$ by
\[
\tilde A v=\left\{
\begin{array}{cl}
A v, & v\in \mathbb{E},\\
0, & v\in \mathbb{E}^\perp.
\end{array}
\right.
\]
Then, it is easy to see that all statements (i1)-(i4), in which we replace the space $\mathbb{E}$ with $\tilde{ \mathbb{E}}$,
the subspaces $V, \Xi_j\subset \mathbb{E}$ with $\tilde V, \tilde \Xi\subset \tilde {\mathbb{E}}$,
the parallelepipeds $\Pi$, $\Pi\cap V$ and their faces $\Pi_j$, $F_j$ with the polyhedra  $\tilde \Pi$, $\tilde P\cap \tilde V$
and their faces $\tilde P_j$, $\tilde F_j$, and the matrix $A$ with the matrix $\tilde A$, remain valid. That is,

\begin{itemize}
	
	\item[(j1)]
	$
	B_j\in \tilde V\cap \overset{\circ}{(\tilde \Pi_1\cap \cdots \cap \tilde \Pi_j)}$ for
	$j=1,\ldots,S$.
	
	\item[(j2)]
	The faces $\tilde F_j=\tilde V\cap \tilde \Pi_1\cap \cdots \cap \tilde \Pi_j $ of $\tilde \Pi\cap \tilde V$ satisfy conditions \eqref{circ}, \eqref{dim}.
	
	\item[(j3)]
	The subspace $\tilde V$ and the hyperplanes
	$
	\tilde \Xi_j$,  $j=1,\ldots,S,
	$
	are in general linear position.
	
	\item[(j4)] The vectors $f_j$ satisfy $f_j={\rm proj}_{\tilde V\cap\tilde \Xi_0\cap \tilde \Xi_1\cap\cdots\cap \tilde\Xi_j} \tilde {A}^\frac12 k_0\ne 0$ for $j=0,\ldots,S-1$.
	
\end{itemize}

Finally, consider the Moreaux process corresponding to the system of $\tilde m $ springs, in which the springs with $(ij)\in \tilde M\setminus M$
have small Young's moduli $a_{ij}$. For this system, we denote by $\hat k_0$, $\hat W, \hat V\subset \tilde{\mathbb{E}}$ and $\hat A$ the counterparts of the vector $k_0$, the subspaces $W, V$ and the matrix $A$ (cf.~Section 2.4). Furthermore, we denote by $\hat \Pi$, $\hat \Pi_j$ the parallelepiped
\[
\hat \Pi=\{v\in \tilde{\mathbb E}: |v_{ij}|\le r_{ij}/\sqrt{a_{ij}},\ (ij)\in \tilde M\}
\]
and its facets
\begin{equation}\label{ad4}
\hat \Pi_j=\{v\in\hat \Pi: v_{1,i_{j-1}+1}=r_{{1,i_{j-1}+1}}/\sqrt{a_{1,i_{j-1}+1}} \}.
\end{equation}
Note that the subspace $\hat V$ is a small perturbation of the subspace $V$
in the sense that the intersections of $V$ and $\hat V$ with any ball centered at the origin
can be made arbitrarily close by making the coefficients $a_{ij}$, $(ij)\in\tilde M\setminus M$,
sufficiently small. Similarly, the sets $\hat \Pi$, $\hat \Pi_j$ are small perturbations of $\Pi$, $\Pi_j$, respectively.
Also, the matrix $\hat A$ and hence the vector $\hat f_0=\hat A^\frac12 \hat k_0$ are small perturbations of the
matrix $\tilde A$ and the vector $f_0=\tilde A^\frac12 k_0$, respectively, because $\hat k_0^{i,j}=k_0^{i,j}$ for $(ij)\in M$.
Therefore, properties (j3), (j4) imply that
the subspace $\hat V$ and the hyperplanes
$
\tilde \Xi_j$,  $j=1,\ldots,S,
$
are in general linear position; each subspace $\hat V\cap\tilde \Xi_0\cap \tilde \Xi_1\cap\cdots\cap \tilde\Xi_j$
is a small perturbation of the subspace $\tilde V\cap\tilde \Xi_0\cap \tilde \Xi_1\cap\cdots\cap \tilde\Xi_j$;
and, each vector $\hat f_j={\rm proj}_{\hat V\cap\tilde \Xi_0\cap \tilde \Xi_1\cap\cdots\cap \tilde\Xi_j} \hat {A}^\frac12 \hat k_0$
is a small perturbation of the vector $\tilde f_j\ne0$, hence $\hat f_j\ne0$. Now, since $B_j$ is the intersection point of the ray
$\{y=B_{j-1} + d f_{j-1}: d>0\}$ with the parallelepiped $\hat\Pi_1\cap\cdots\cap \hat\Pi_j$ for every $j=1,\ldots,S$, property (j1)
implies that the relation
\[
\psi_j(x):= \{y=x + d \hat f_{j-1}: d>0\}\cap (\hat\Pi_1\cap\cdots\cap \hat\Pi_j)
\]
uniquely defines a point $\psi_j(x)\in \overset{\circ}{(\hat\Pi_1\cap\cdots\cap \hat\Pi_j)}$
whenever a point $x\in \hat\Pi_1\cap\cdots\cap \hat\Pi_{j-1}$ is sufficiently close to $B_{j-1}$
%here $x$ is simply close to $B_0=0$ for $j=1$
and the vector $\hat f_{j-1}$ is sufficiently close to $f_{j-1}$. Therefore,
for sufficiently small $a_{ij}$, $(ij)\in\tilde M\setminus M$, the points
$\hat B_0=0, \hat B_1=\psi_1(\hat B_0), \hat B_2=\psi_2(\hat B_1),\ldots, \hat B_S=\psi_S(\hat B_{S-1})$ satisfy
\begin{equation}\label{add3}
\hat B_j\in \overset{\circ}{(\hat\Pi_1\cap\cdots\cap \hat\Pi_j)},\qquad \hat B_j=\hat B_{j-1}+(\hat d_{j}-\hat d_{j-1})\hat f_{j-1}, \ \ \ j=1,\ldots,S,
\end{equation}
with $0=\hat d_0<\hat d_1<\ldots<\hat d_S$, and the polyline $\hat \gamma=\hat B_0\hat B_1\cdots\hat B_S$ is a small perturbation of the polyline
$\gamma=B_0B_1\ldots B_S$. Since $\hat f_j\in \hat V$ for all $j$, we also see that the vertex $\hat B_j$ of $\hat\gamma$ belongs to the open face
$\overset{\circ}{\hat F}_j=\hat V\cap \overset{\circ}{(\hat\Pi_1\cap\cdots\cap \hat\Pi_j)}$ of the polytope $\hat V\cap\hat \Pi$
and the faces $\hat F_j$ satisfy conditions \eqref{circ}, \eqref{dim}.

It remains to show that the set
\[
\hat  \Omega = \{ y = \tau_1 \hat B_0\hat B_1 +\cdots + \tau_S \hat B_{S-1}\hat B_S, \ |\tau_i|\le 1, \ i=1,\ldots,S\}\subset \hat V
\]
(cf.~\eqref{omega}) belongs to $\hat \Pi$. To this end, denote by $e_{ij}$ the standard basis of vectors $e_{i,j}$ in $\tilde {\mathbb {E}}$:
\[
e_{i,j}^{i',j'}=\left\{ \begin{array}{cl}
1,& (i'j')=(ij),\\
0, & {\rm otherwise.}
\end{array}\right.
\]
For $(i,j)=(1,i_n+1)$, $n=0,\ldots,S-1$, consider the estimate
\[
|\langle e_{1,i_n+1}, \sum_{j=1}^S \tau_j \hat B_{j-1}\hat B_j \rangle|\le \sum_{j=1}^S |\langle e_{1,i_n+1}, \hat B_{j-1}\hat B_j \rangle|
\]
with $|\tau_j|\le 1$.
Here $\hat B_{n+1},\ldots \hat B_S\in \hat \Pi_{n+1}$, hence $\langle e_{1,i_n+1}, B_{j-1}B_j\rangle =0$ for $j=n+2,\ldots S$
and therefore
\begin{equation}\label{ta}
|\langle e_{1,i_n+1}, \sum_{j=1}^S \tau_j \hat B_{j-1}\hat B_j \rangle|\le \sum_{j=1}^{n+1} |\langle e_{1,i_n+1}, \hat B_{j-1}\hat B_j \rangle|.
\end{equation}
But relations \eqref{add1}, \eqref{add2} imply that
$\langle e_{1,i_n+1},f_{j-1}\rangle >0$ for $1\le j\le n+1\le S$
and since $\hat f_{j-1}$ is a small perturbation of the vector $f_{j-1}$,
\[
\langle e_{1,i_n+1},  \hat f_{j-1} \rangle >0, \qquad 1\le j\le n+1\le S.
\]
Combining this estimate with the second relation in \eqref{add3} results in the relation
$\langle e_{1,i_n+1},  \hat B_{j-1}\hat B_j \rangle >0$ for $j=1,\ldots,n+1$, hence \eqref{ta} is equivalent to
\begin{equation}\label{comb1}
|\langle e_{1,i_n+1}, \sum_{j=1}^S \tau_j \hat B_{j-1}\hat B_j \rangle|\le \sum_{j=1}^{n+1} \langle e_{1,i_n+1}, \hat B_{j-1}\hat B_j \rangle= \langle e_{1,i_n+1}, \hat B_{0}\hat B_{n+1} \rangle
=\frac{r_{1,i_n+1}}{\sqrt{a_{1,i_n+1}}},
\end{equation}
where the last equality follows from $B_{n+1}\in \hat \Pi_{n+1}$ and \eqref{ad4}.

Now, consider $(ij)\ne (1,i_{n+1})\in M$, i.e.~$(i,j)=(i,{q_n(i)})$ for some $0\le n\le S-1$, $ i_n+1\le i\le i_{n+1}$.
Formulas \eqref{add1}, \eqref{add2} imply that $\langle e_{ii'}, f_{j-1}\rangle\ge 0$ for all $(ii')\in M$, $j=1,\ldots, S$
and $\langle e_{i,q_n(i)},f_{j-1}\rangle =0$ for $j>n+1$.
Therefore, taking into account that $B_{j-1}B_j=(d_j-d_{j-1})f_{j-1}$ and arguing exactly in the same way as above,
we see that the unperturbed polyline $\gamma$ satisfies
\[
|\langle e_{i,q_n(i)}, \sum_{j=1}^S \tau_j B_{j-1} B_j \rangle|\le \sum_{j=1}^{n+1} \langle e_{i,q_n(i)}, B_{j-1} B_j \rangle. %= \langle e_{i,q_n(i)}, B_{0} B_{n+1} \rangle.
\]
Furthermore, formulas \eqref{add1}, \eqref{add2} and $B_{j-1}B_j=(d_j-d_{j-1})f_{j-1}$ imply
\[
{\langle e_{i,q_n(i)},B_{j-1}B_{j}\rangle}{\sqrt{a_{i,q_n(i)}}}={\langle e_{1,i_n+1},B_{j-1}B_{j}\rangle}{\sqrt{a_{1,i_n+1}}}, \qquad j=1,\ldots,n+1,
%(d_{n+1}-d_n)\langle e_{i,q_n(i)},f_{n}\rangle
\]
hence
\[
|\langle e_{i,q_n(i)}, \sum_{j=1}^S \tau_j B_{j-1} B_j \rangle|\le %\sqrt{\frac{a_{1,i_n+1}}{a_{i,q_n(i)}}}\sum_{j=1}^{n+1} \langle e_{1,i_n+1}, B_{j-1} B_j \rangle=
\langle e_{1,i_n+1}, B_{0} B_{n+1} \rangle \sqrt{\frac{a_{1,i_n+1}}{a_{i,q_n(i)}}},
\]
where, due to $B_{n+1}\in \hat \Pi_{n+1}$,
\[
{\langle e_{1,i_n+1},B_{0}B_{n+1}\rangle}{\sqrt{a_{1,i_n+1}}}=r_{{1,i_{n}+1}}.
\]
Therefore, \eqref{zz} implies
\[
|\langle e_{i,q_n(i)}, \sum_{j=1}^S \tau_j B_{j-1} B_j \rangle|<\frac{r_{i,q_n(i)}}{\sqrt{a_{i,q_n(i)}}}.
\]
Since $\hat \gamma$ is a small perturbation of $\gamma$, a similar estimate is true for $\hat\gamma$:
\begin{equation}\label{comb2}
|\langle e_{i,q_n(i)}, \sum_{j=1}^S \tau_j \hat B_{j-1} \hat B_j \rangle|<\frac{r_{i,q_n(i)}}{\sqrt{a_{i,q_n(i)}}},\qquad i_n+1\le i\le i_{n+1}.
\end{equation}

Finally, for $ii')\in\tilde M\setminus M$,
\begin{equation}\label{comb3}
|\langle e_{i,i'}, \sum_{j=1}^S \tau_j \hat B_{j-1} \hat B_j \rangle|<\frac{r_{i,i'}}{\sqrt{a_{i,i'}}},
\end{equation}
because $a_{i,i'}$ is small.
Combining \eqref{comb1}, \eqref{comb2} and \eqref{comb3}, we obtain $\hat \Omega\subset \hat \Pi$.
This completes the proof of Theorem 2.

%\bigskip
%{\bf Aside.} Additional formulas:
%
%$\overset{\circ}{E}_j$ consists of vectors $v$ satisfying \eqref{vvvv} and, additionally,
%$v_{i,q_n(i)}=0$ for $i\in \{1, i_n+1,i_n+2,\ldots, i_{n+1}\}, \ n< j$.
%
%$n_j$ is defined by
%\[
%n_j^{i,q_{j-1}(i)} \sqrt{a_{i,q_{j-1}(i)}} =
%\left\{
%\begin{array}{cl}
%1, & i\in \{1, i_{j-1}+1,i_{j-1}+2,\ldots, i_{j}\},\\
%0, & {\rm otherwise}.
%\end{array}
%\right.
%\]
%
%$c_j$ is defined by
%\[
%c_{n+1}\left(a_{1,i_n+1}+\sum_{i=i_n+1}^{i_{n+1}-1}a_{i,i+1}+a_{i_{n+1},N}\right)=\tilde r_{n},\qquad\quad n=0,\ldots,S-1.
%\]

%\subsection{Queueing theory interpretation}
%Here we briefly discuss another interpretation of Fig.~\ref{fig1}.
%We consider a stationary flow of customers from node 1 to node 5 through
%a network of channels represented by the links between the nodes.
%The flux of customers through the link $(ij)$ with $i<j$ is $\sigma_{ij}$ (a counterpart of stress in the mechanical model).
%As the flow is stationary, the sum of signed fluxes $(-1)^k_{ij} \sigma_{ij}$ over all the links starting at a node $i$
%(where $k_{ij}=0$ for $j>i$ and $k_{ij}=1$ for $j<i$) is zero for each node $i=2,3,4$.

\section{Conclusions}
We considered a set of Prandtl's elastic-ideal plastic elements, which are arranged into a network and deform quasistatically
when the distance between two nodes is varied according to a given law (input). As in the general setting of the model of Moreau,
no {\em a priori} constraint was imposed on the topology of the network. We defined the loading curve $\phi$ for the corresponding sweeping process
as a graph of the solution corresponding to the zero initial condition and an increasing input. It was shown that if this curve satisfies simple geometric conditions,
then the structure of hysteresis loops of the model in the space of stresses extended by one dimension representing the input
is similar to the structure of hysteresis loops of the PI model. Furthermore, the relationship between the input and the varying stress
of each Prandtl's element in the network is given by a PI operator ${\mathcal P}_{\phi_i}$,
and the loading curve $\phi_i$ of this operator is the corresponding projection of the loading curve $\phi$ of the sweeping process.
The question of how the geometric conditions of the main theorem can be related to the topology and parameters of the network
in general remains open. However, we showed that these conditions are satisfied for any Moreau network obtained by a small perturbation of a PI model.
In other words, they are satisfied for any network with sufficiently small coupling. This is in line with the results from \cite{xx3} obtained by a different approach based on 
the composition property of PI operators.

In this work, we considered networks of springs aligned along a straight line and used the distance between two nodes A and B as the input.
However, the results can be easily generalized to networks of Prandtl's elements
connecting nodes in a two- or three-dimensional space. It would be interesting to consider other types of scalar-valued inputs,
for example, the force applied at the node A. In this case, the input $Z(t)$ of the sweeping process 
is a set with changing size and shape and hence is more complicated than the inputs we considered.
This will be the subject of future work.

%Queueing theory.

\section*{Acknowledgments}
The author thanks Ivan Gudoshnikov for multiple stimulating discussions
of the Moreau model.

The author acknowledges the support of NSF through grant DMS-1413223.


\begin{thebibliography}{99}

\bibitem{mroz}
M Brokate, P Krejci, D Rachinskii,
Some analytical properties of the multidimensional continuous Mroz model of plasticity,
Control and Cybernetics 27 (2), 1998, 199-215
	
	\bibitem{xx2}
	A Krasnosel'skii, D Rachinskii, 
	On a bifurcation governed by hysteresis nonlinearity,
	Nonlinear Differential Equations and Applications NoDEA 9 (1), 2002, 93-115
	
	\bibitem{xx3}
	P Krejci, H Lamba, S Melnik, D Rachinskii,
	Analytical solution for a class of network dynamics with mechanical and financial applications, 
	Physical Review E 90 (3), 2014, 032822
	
	\bibitem{xx4}
	P Krejci, H Lamba, S Melnik, D Rachinskii,
	Kurzweil integral representation of interacting Prandtl-Ishlinskii operators,
	Discrete Contin. Dyn. Syst., Ser. B 20, 2015, 2949-2965
	
	\bibitem{xx5}
	P Krejci, H Lamba, GA Monteiro, D Rachinskii,
	The Kurzweil integral in financial market modeling,
	Mathematica Bohemica 141 (2), 2016, 261-286
	
	\bibitem{KP}
	M Krasnosel'skii, A. Pokrovskii, Systems with Hysteresis, Springer, 1989
	
	\bibitem{bs}
	M. Brokate, J. Sprekels, Hysteresis and Phase Transitions, Springer, 1996
		
	\bibitem{xx1}
	P Krejci, J P O'Kane, A Pokrovskii, D Rachinskii, 
	Properties of solutions to a class of differential models incorporating Preisach hysteresis operator,
	Physica D: Nonlinear Phenomena 241 (22), 2012, 2010-2028
	
	\bibitem{xx6}
	D Rachinskii, 
	Realization of arbitrary hysteresis by a low-dimensional gradient flow,
	Discrete Contin. Dyn. Syst., Ser. B 21 (1), 2016, 227-243
	
\end{thebibliography}
\end{document}